\documentclass[article]{IEEEtran}
\usepackage{epsfig}
\usepackage{graphicx}
\usepackage{psfrag}
\usepackage{ifpdf}
\usepackage{cite}
\usepackage{siunitx}
\usepackage[utf8x]{inputenc}
\usepackage[top=1in, bottom=1in]{geometry}
\usepackage{fancyhdr}
\usepackage{lastpage}
\usepackage{oldstyle}
\usepackage{textcomp}
\usepackage{color}
\usepackage{placeins}
\usepackage{float}
\usepackage{tabularx,colortbl}
\graphicspath{{figures/Appendix/}}
\usepackage{amssymb} 
\usepackage{amsmath} 
\usepackage{subfigure}
\usepackage{algorithm}
\usepackage{algpseudocode}
\usepackage{listings}
\definecolor{mygreen}{RGB}{28,172,0} 
\definecolor{mylilas}{RGB}{170,55,241}

\usepackage{array}
\usepackage{setspace}

\usepackage{amsmath,graphicx}
\usepackage{epstopdf}
\usepackage{morefloats}
\usepackage[square, comma, sort&compress, numbers]{natbib}

\onecolumn

\begin{document}

\begin{center}
 \LARGE Denosing Using Wavelets and Projections onto the $\ell_1$-Ball\\[10pt]
 \normalsize\today
\end{center}

\begin{center}
\normalsize A. Enis Cetin, M. Tofighi\\
Dept. of Electrical and Electronic Engineering, Bilkent University,  Ankara, Turkey\\
cetin@bilkent.edu.tr, tofighi@ee.bilkent.edu.tr\\
\end{center}

\vspace{0.5cm}


Both wavelet denoising and denosing methods using the concept of sparsity are based on soft-thresholding. In sparsity based denoising methods, it is assumed that the original signal is sparse in some transform domains such as the wavelet domain and the wavelet subsignals of the noisy signal are projected onto $\ell_{1}$-balls to reduce noise. In this lecture note, it is shown that the size of the $\ell_{1}$-ball or equivalently the soft threshold value can be determined using linear algebra. The key step is an orthogonal projection onto the epigraph set of the $\ell_{1}$ norm cost function.

\IEEEpeerreviewmaketitle

\IEEEpubidadjcol
In standard wavelet denoising, a signal corrupted by additive noise is wavelet transformed and resulting wavelet subsignals are soft and/or hard thresholded. After this step the denoised signal is reconstructed from the thresholded wavelet subsignals \cite{mallat,Donoho}. Thresholding the wavelet coefficients intuitively makes sense because wavelet subsignals obtained from an orthogonal or biorgthogonal wavelet filterbank exhibit large amplitude coefficients only around edges or change locations of the original signal. Other small amplitude coefficients should be due to noise. Many other related wavelet denoising methods are developed based on Donoho and Johnstone's idea, see e.g. \cite{Donoho,Mina,Rousseeuw,mallat,vetterli,Johnstone}. Most denoising methods take advantage of sparse nature of practical signals in wavelet domain to reduce the noise \cite{PesquetICASSP,GlobalSIP2013,Kose201463,Pesquet,Duchi,Bar07}.  

Consider the following basic denoising framework. Let $v[n]$ be a discrete-time signal and $x[n]$ be a noisy  version of $v[n]$:
\begin{equation}
x[n]=v[n] + \xi[n], \ \ \ n=0,1,2,\dots,N-1.
\end{equation}
where $\xi[n]$ is the additive, i.i.d, zero-mean, white Gaussian noise with variance  $\sigma^2$. An L-level discrete wavelet transform of $x[n]/ \sqrt{N}$ is computed and the lowband signal $x_L$ and wavelet subsignals $w_{1},w_{2},\dots,w_{L}$ are obtained as shown in Fig. 1. After this step, wavelet subsignals are soft-thresholded as shown in Fig. \ref{wavcoeff}. The soft threshold, $\theta_{i}$, can be selected in many ways. One possible choice is 
\begin{equation}
\theta_{i} =\gamma.\sigma.\sqrt{2log(N)/N},
\end{equation}
 where $\gamma$ is a constant \cite{Donoho}. The problem with this threshold is that the noise variance $\sigma^2$ has to be known or properly estimated from the observations, $x[n]$. 

Another way to denoise the wavelet subsignals $w_{1},w_{2},\dots,w_{L}$ is to project them onto $\ell_1$-balls. As pointed out above denoising is possible with the assumption that wavelet subsignals are also sparse signals. Projection $w_{pi}[n]$ of $w_{i}[n]$ onto an $\ell_1$-ball is obtained as follows:
\begin{align}
\label{eq:eq7}
w_{pi}&= \text{argmin}\|w_{i}-w\|_{2} ^2 \\ \nonumber
& \text{such that}  \sum_n |w[n]| \leq d_{i},
\end{align}
where $d_{i}$ is the size of the $i$-th $\ell_{1}$-ball. This minimization problem was studied by many researchers and computationally efficient algorithms were developed (see e.g., \cite{Duchi}). The projection vector $w_{pi}$ is basically obtained by soft-thresholding as in ordinary wavelet denoising. After orthogonal projection onto an $\ell_{1}$-ball each wavelet coefficient is modified as follows:
\begin{equation}
\label{eq:eq6}
w_{pi}[n]= \text{sign}(w_{i}[n]).\text{max}\big\{|w_{i}[n]|- \theta_{i},0\big\},
\end{equation}
where $\theta_{i}$ is a constant whose value is determined according to the size of the $\ell_{1}$-ball, $d_{i}$, as described in Algorithm 1. Equation \ref{eq:eq6} is basically soft-thresholding with $\theta_{i}$ as the threshold value. Other computationally efficient algorithms capable of computing the projection  vector $w_{pi}$ in $O(K)$ time are described in \cite{Duchi}.
Projection operations onto $\ell_1$-ball will force small valued wavelet coefficients to zero and retain the edges and sharp variation regions of the signal because wavelet subsignals have large amplitudes corresponding to edges in most natural signals. As in standard wavelet denoising methods the low-band subsignal $x_L$ is not processed because $x_L$ is not a sparse signal for most practical signals. 

In standard wavelet denoising, noise variance has to be estimated to determine the soft-threshold value. In this case, the size of the $\ell_1$-ball $d_i$ in (\ref{eq:eq7}) or equivalently $\theta_i$ in (\ref{eq:eq6}) has to be estimated. Another parameter that has to be determined in both standard wavelet denoising and the $\ell_1$-ball based denoising is the number of wavelet decomposition levels. In the next two sections we describe how the size of the $\ell_1$-ball and the number of wavelet decomposition levels can be determined.

\vspace{-0.03cm}

\section{Estimation of Denoising Thresholds Using the Epigraph Set of $\ell_1$-ball}
The soft-threshold $\theta_i$ is related with the size of $\ell_1$-ball as described in Algorithm 1. The size of the $\ell_1$-ball can vary between 0 and $d_{max,i}$ which is determined by the boundary of the $\ell_1$-ball going through the wavelet subsignal $w_{i}[n]$:
\begin{align}
d_{max,i}=\sum_n sign(w_{i}[n]) w_{i}[n],
\end{align}
where $sign(w_{i}[n])$ is the sign of $w_{i}[n]$. Orthogonal projection of $w_{i}[n]$ onto a ball with $d=0$ produces an all-zero result. On the other hand, projection of $w_{i}[n]$ onto a ball with size $d_{max,i}$, does not change $w_{i}[n]$ because $w_{i}[n]$ is on the boundary of the $\ell_1$-ball. Therefore,  the ball size $z$ must satisfy the inequality $0< z <d_{max}$, for denoising. This $\ell_1$-ball condition can be expressed  as follows: 
\begin{align}
g(\textbf{w}) = \sum_{k=0}^{K-1}|w[k]|\leq z,
 \label{1}
\end{align}
where $K$ is the length of the wavelet subsignal $w_{i}[n]$ and ${\textbf{w}} = [w[0], w[1], \dots, w[K-1]]$. This condition corresponds to the epigraph set of the $\ell_1$-ball in $\mathbb{R}^{K+1}$ \cite{PesquetICASSP, Pesquet}. In (\ref{1}) there are $\mathbb{R}^{K+1}$ variables. These are $w_{i}[0], \dots, w_{i}[K-1]$, and $z$. The epigraph set of $\ell_1$-cost function {${g(\textbf{w})\leq z}$} in $\mathbb{R}^3$ is shown in Fig. \ref{2D_L1}.

\algnewcommand{\Inputs}[1]{%
  \State \textbf{Inputs:}
  \Statex \hspace*{\algorithmicindent}\parbox[t]{.8\linewidth}{\raggedright #1}
}
\algnewcommand{\Initialize}[1]{%
  \State \textbf{Initialize:}
  \Statex \hspace*{\algorithmicindent}\parbox[t]{.8\linewidth}{\raggedright #1}
}
\algnewcommand{\Output}[1]{%
  \State \textbf{Output:}
  \Statex \hspace*{\algorithmicindent}\parbox[t]{.8\linewidth}{\raggedright #1}
}

\begin{algorithm}[!htb]
  \caption{Order ($Klog(K)$) algorithm implementing projection onto
    the $\ell_1$-ball with size $d_i$.}
  {\fontsize{9}{5}
  \begin{algorithmic}[1]
 \Inputs{A vector $w_{i}[n],~n=0, 1, \dots, K-1$ and scalar $d_i>0$}
     \Initialize{Sort the entries of $|w_{i}[n]|$ and obtain the rank ordered sequence $\mu_{1} \geq \mu_{2}\geq, \dots,\geq \mu_{K}$}
\begin{align}
\theta_{i}&= \frac{1}{\rho} \big(\sum_{n=1}^\rho \mu_n-d_{i} \big) \  \\ \nonumber
&\text{such that} ~ ~ \rho=max \big\{j\in \{0, 1, 2, \dots, K-1\}:  \\ \nonumber 
&\mu_{j} - \frac{1}{j}\big(\sum_{r=1}^{j}\mu_{r}-d_i \big)>0  \big\} 
\end{align}
 \Output{$w_{pi}[n]~=~sign(w_{i}[n]).max \{|w_{i}[n]|~-~\theta_{i}, 0\},$
 $ n = 0, 1, 2, \dots, K-1$}
  \end{algorithmic} }
\end{algorithm}

 \begin{figure}[!htb]
 \centering
 \includegraphics[scale=0.15]{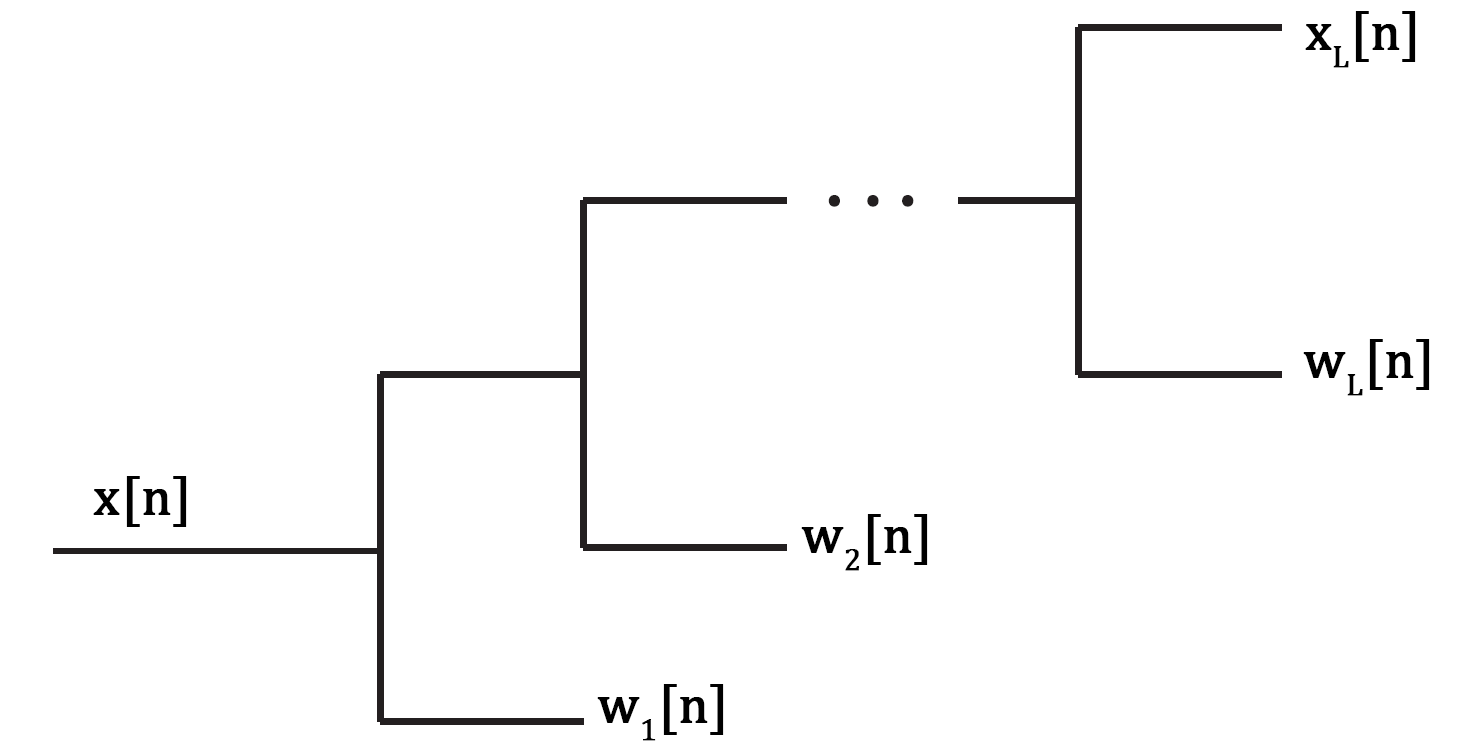}
 \caption{L-level dyadic wavelet decomposition of the signal x.}
 \label{sema1}
 \end{figure}
By orthogonal projecting the wavelet subsignal $[w_{i}[n], 0]^{T} := [w_{i}[0], \dots, w_{i}[K-1], 0]^{T}$ onto the epigraph set it is possible to determine all of the $\mathbb{R}^{K+1}$ unknowns, $w_{pi}[n],~n = 0, 1, \dots, K-1$, and $z_p$, as graphically illustrated in Fig. \ref{L1norm}. The projection vector $[w_{pi}[n], d]^T$ is unique and the closest vector to the wavelet subsignal $[w_{i}[n], 0]^T$ in $\mathbb{R}^{K+1}$ because the epigraph set is a closed and convex set. The projection onto the epigraph set can be computed in two steps. In the first step, $[w_{i}[n], 0]^{T}$ is projected onto the nearest boundary hyperplane of the epigraph set which is 
\begin{align}
\sum_{n=0}^{K-1}sign(w_{i}[n]).w[n]-z=0.
\end{align}
 \begin{figure}[!htb]
 \centering
 \includegraphics[scale=0.2]{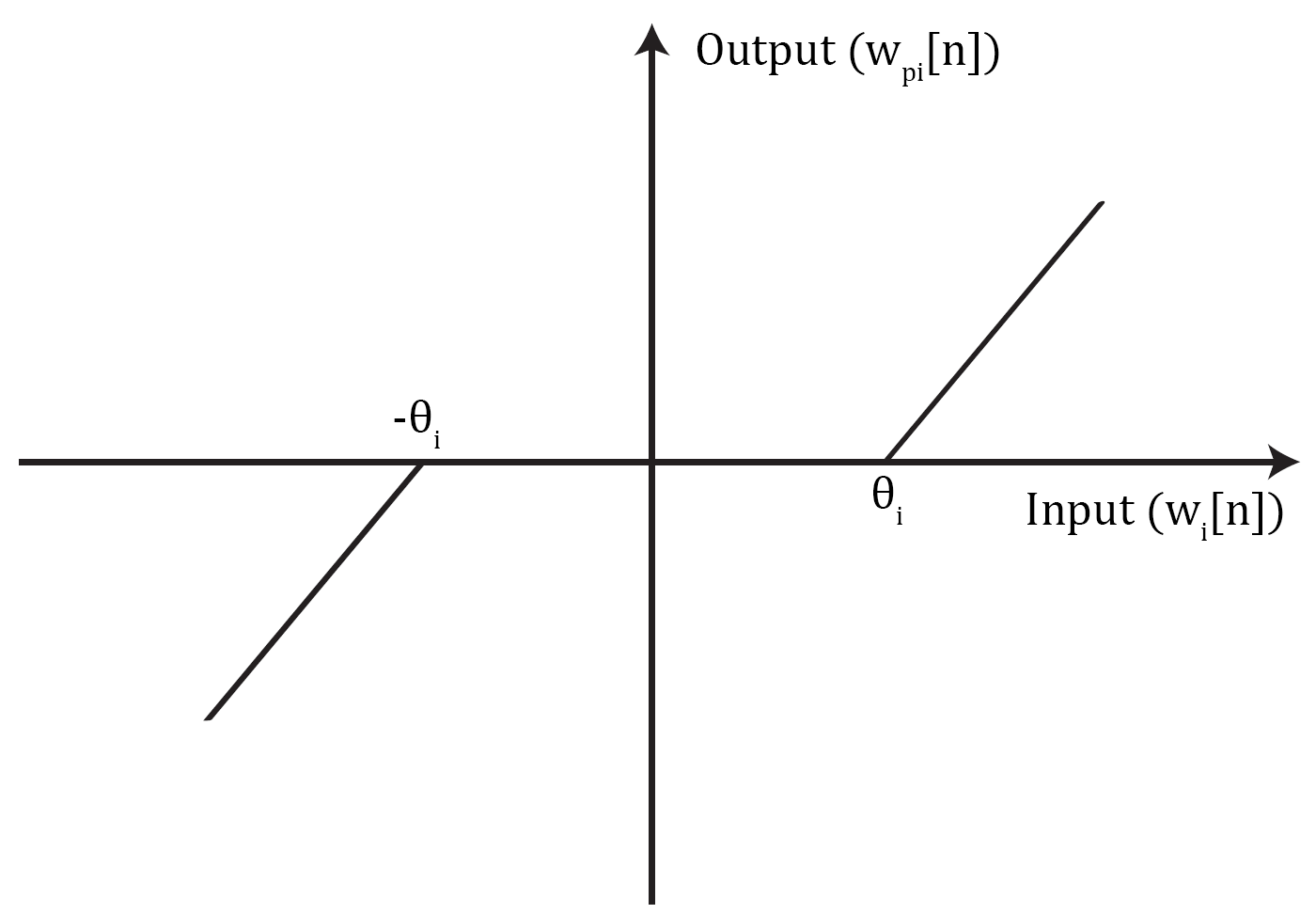}
 \caption{Soft-thresholding of wavelet coefficients.}
 \label{wavcoeff}
 \end{figure}
The projection signal $w_{pi}[n]$ onto the hyperplane in $\mathbb{R}^{K+1}$ is determined as follows:
\begin{equation}
\begin{split}
w_{pi}[n]= w_{i}[n] + \frac{0.z - \sum_{n=0}^{K-1}sign(w_{i}[n])w_{i}[n]}{K+1} sign(w_{i}[n]) \\
 \quad n = 0, 1, \dots, K-1,
\end{split}
\end{equation}
and
\begin{equation}
\label{eq:eq9}
z_p = 0 + \frac{\sum_{n=0}^{K-1}sign(w_{i}[n])w_{i}[n]}{K+1}.
\end{equation}
This orthogonal projection also determines the size of the ball:
\begin{equation}
\label{eq:eq10}
d_i = \sum_{n=0}^{K-1}sign(w_{i}[n])w_{pi}[n],
\end{equation}
because the projection vector $w_{pi}[n]$, $n = 0, 1, \dots, K-1$ must be on the $K$-dimensional hyperplane with weights $sign(w_{i}[n])$. This is graphically illustrated in Fig. \ref{2D_L1} (view from the top).
 \begin{figure}[!htb]
 \centering
 \includegraphics[scale=0.22]{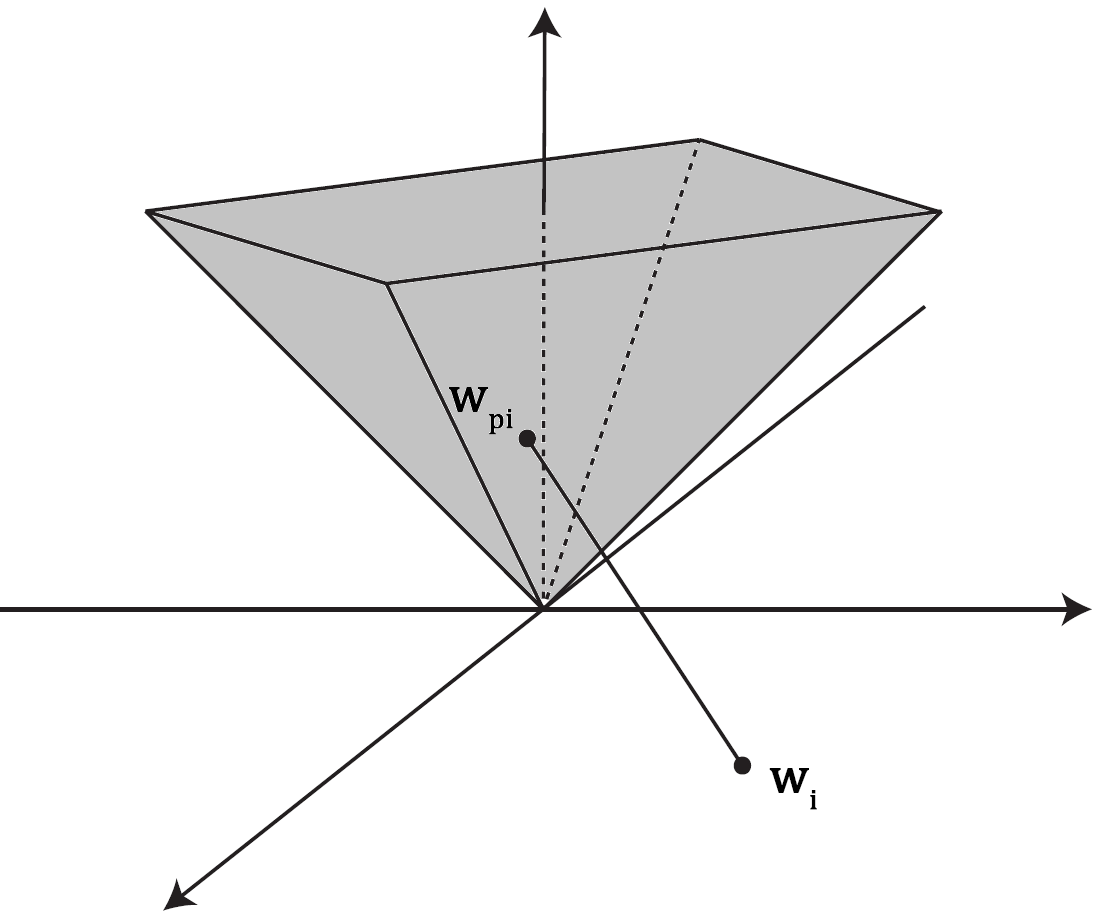}
 \caption{Projection of $w_{i}[n]$ onto the epigraph set of $\ell_1$-norm cost function: $z \geq \sum_{n=0}^{K-1}|w[k]|$}
 \label{L1norm}
 \end{figure}
The projection $w_{pi}[n]$ may or may not be on the epigraph set of $\ell_1$-ball. If the signs of the projection signal $w_{pi}[n]$ entries are the same as $w_{i}[n]$ for all $n$ then the $w_{pi}[n]$ is on the epigraph set, otherwise $w_{pi}[n]$ is not on the $\ell_1$-ball as shown in Fig. \ref{2D_L1}. If $w_{pi}[n]$ is not on the $\ell_1$-ball we can still project $w_{i}[n]$ onto the $\ell_1$-ball using Algorithm 1 or Duchi et al's $\ell_1$-ball projection algorithm \cite{Duchi} using the values of $d_i$ determined in Eq. (\ref{eq:eq10}). This constitutes the second step of epigraph projection operation. 

In summary, we have the following two steps: (i) Project $w_{i}[n]$ onto the boundary hyperplane and determine $d_i$. (ii) If $sign(w_{i}[n]) = sign(w_{pi}[n])$ for all $n$, $w_{pi}[n]$ is the projection vector. Otherwise use $d_i$ value in Algorithm 1 to determine the final projection vector. Vector $w_{pi}$ is the projection of $w_{i}$ onto the $\ell_{1}$-ball, but $w_{pi}$ is only the projection of $w_1$ onto one of the boundary hyperplanes of $\ell_{1}$-ball.
\begin{figure}[!htb]
\centering
\includegraphics[scale=0.22]{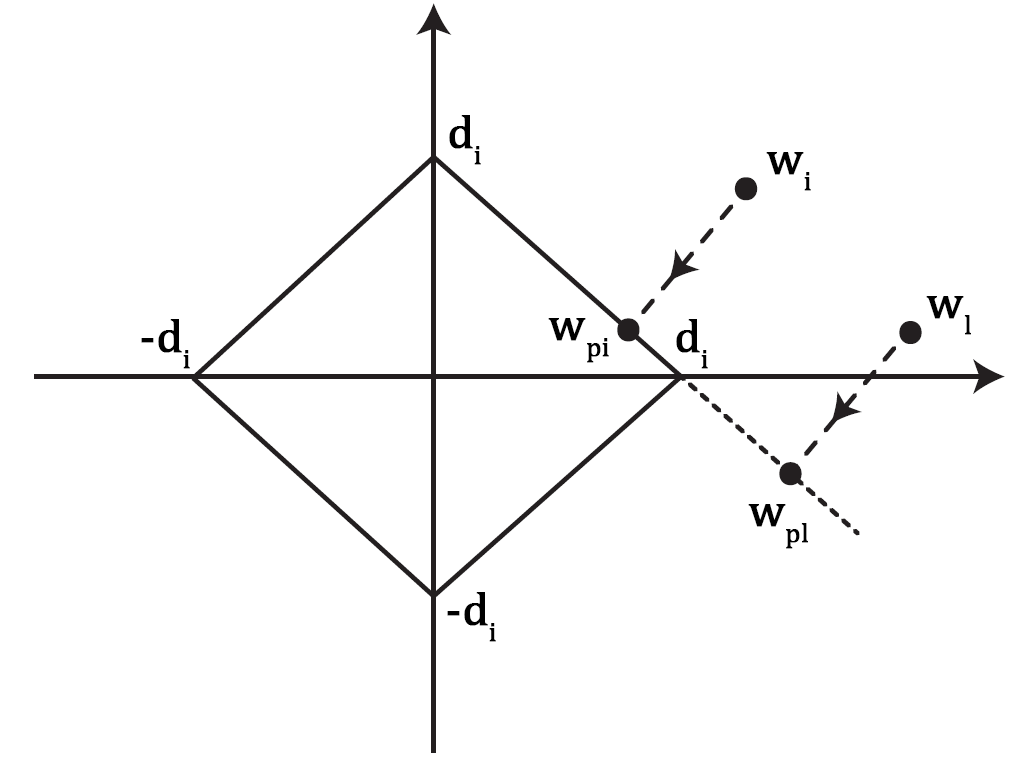}
\caption{Orthogonal projection operation onto a bounding hyperplane of $\ell_1$-ball.}
\label{2D_L1}
\vspace{-0.3cm}
\end{figure}
\section{How to Determine The Number of Wavelet Decomposition Levels}
\label{sec:Wavelet Decomposition Levels}
It is possible to use the Fourier transform of the noisy signal to estimate the bandwidth of the signal. Once the bandwidth $\omega_0$ of the original signal is approximately determined  it can be used to estimate the number of wavelet transform levels and the bandwidth of the low-band signal $x_L$. In an $L$-level wavelet decomposition the low-band signal $x_L$ approximately comes from the $[0, \frac{\pi}{2^L}]$ frequency band of the signal $x[n]$. Therefore, $\frac{\pi}{2^L}$ must be greater than $\omega_0$ so that the actual signal components are not soft-thresholded. Only wavelet subsignals $w_{L}[n], w_{L-1}[n], \dots, w_{1}[n]$, which come from frequency bands  $[\frac{\pi}{2^L}, \frac{\pi}{2^{L-1}}]$, $[\frac{\pi}{2^{L-1}}, \frac{\pi}{2^{L-2}}]$, \dots, $[\frac{\pi}{2}, \pi]$, respectively, should be soft-thresholded in denoising.
\begin{figure*}[!htb]
\centering
\includegraphics[scale=0.15]{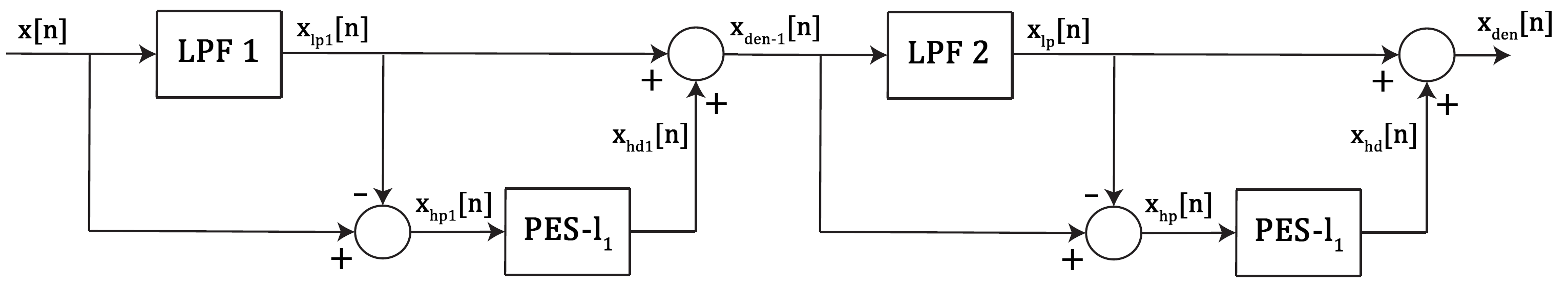}
\caption{Pyramidal filtering based denoising. the high-pass filtered signal is projected onto the epigraph set of $\ell_1$.}
\label{Pyramid}
\end{figure*}
For example, in Fig. \ref{WavLevel}, the magnitude of Fourier transform of $x[n]$ is shown for ``piece-regular" signal defined in MATLAB. This signal is corrupted by zero-mean white Gaussian noise with $\sigma = 10, 20$, and $30 \%$ of the maximum amplitude of the original signal, respectively. For this signal an $L=3$ level wavelet decomposition is suitable because Fourier transform magnitude approaches to the noise floor level after $\omega_{0} = \frac{58\pi}{512}$. It is also a good practice to allow a margin for signal harmonics. Therefore, ($\frac{\pi}{2^{3}}>\frac{58\pi}{512}$) is selected as the number of wavelet decomposition levels.
\begin{figure}[!htb]
\centering
\includegraphics[scale=0.175]{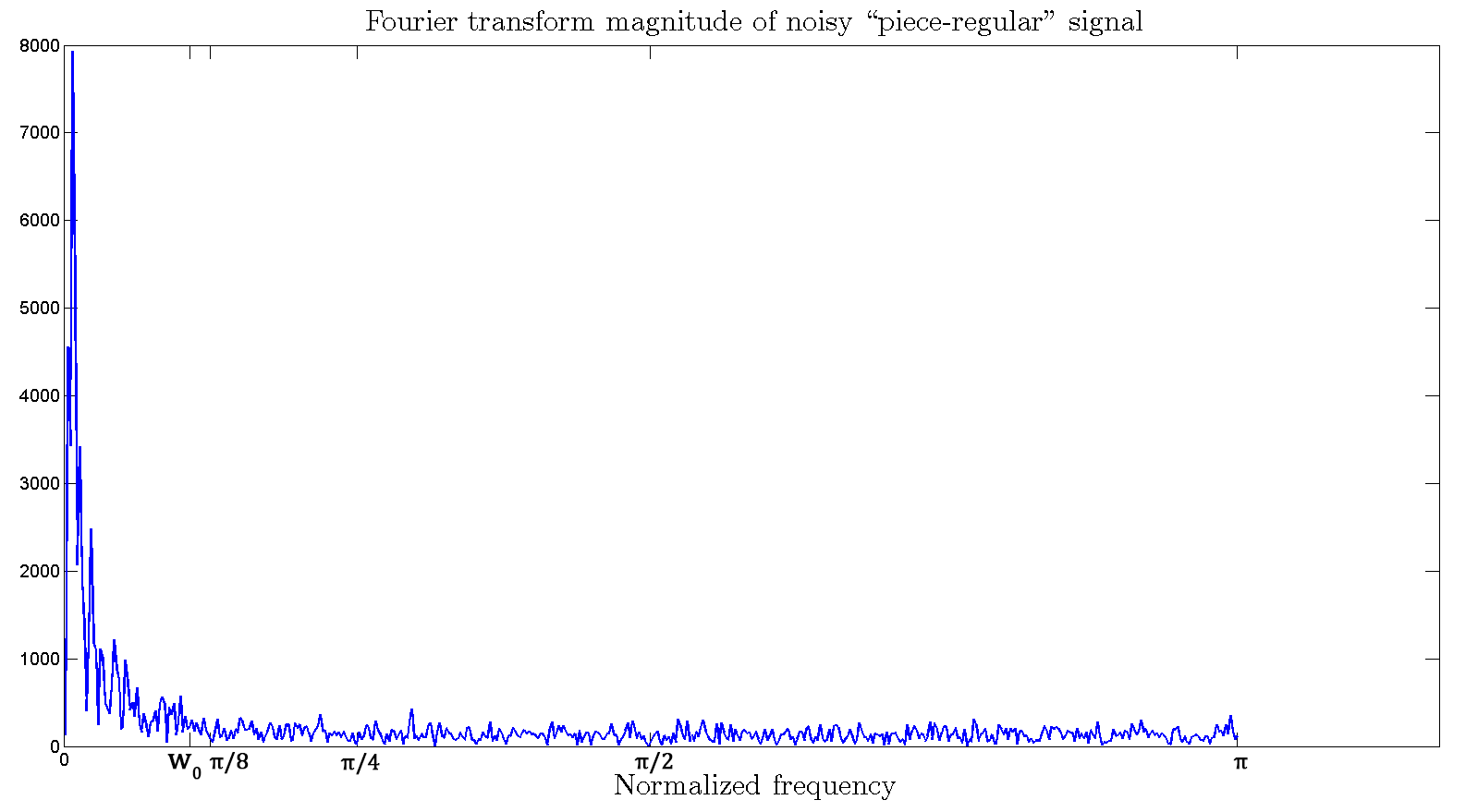}
\caption{Discrete-time Fourier transform magnitude of ``piece-regular" signal corrupted by noise. The wavelet decomposition level L is selected as 3 to satisfy $\frac{\pi}{2^3}>\omega_0$, which is the approximate bandwidth of the signal.}
\label{WavLevel}
\end{figure}
It is also possible to use a pyramidal structure for signal decomposition instead of the wavelet transform. The noisy signal is low-pass filtered with cut-off frequency $\frac{\pi}{8}$ for ``piece-regular" signal and the output $x_{lp}[n]$ is subtracted from the noisy signal $x[n]$ to obtain the high-pass signal $x_{hp}[n]$ as shown in Fig. \ref{Pyramid}. The signal is projected onto the epigraph of $\ell_1$-ball and $x_{hd}[n]$ is obtained. Projection onto the Epigraph Set of $\ell_1$-ball (PES-$\ell_1$), removes the noise by soft-thresholding. The denoised signal $x_{den}[n]$ is reconstructed by adding $x_{hd}[n]$ and $x_{lp}[n]$ as shown in Fig. \ref{Pyramid}. It is possible to use different thresholds for different subbands as in wavelet transform, using a multisatge pyramid as shown in Fig. \ref{Pyramid}. In the first stage a low-pass filter with cut-off $\frac{\pi}{2}$ can be used and $x_{hp1}[n]$ is projected onto the epigraph set of $\ell_1$-ball producing a threshold for the subband $[\frac{\pi}{2}, \pi]$. In the second stage, another low-pass filter with cut-off $\frac{\pi}{4}$ can be used and $x_{hp}[n]$ is projected onto the epigraph set producing a threshold for $[\frac{\pi}{4}, \frac{\pi}{2}]$, etc.

\begin{figure*}[htb!]
\centering
\subfigure[Original signal]
{\label{original}\includegraphics[width=75mm]{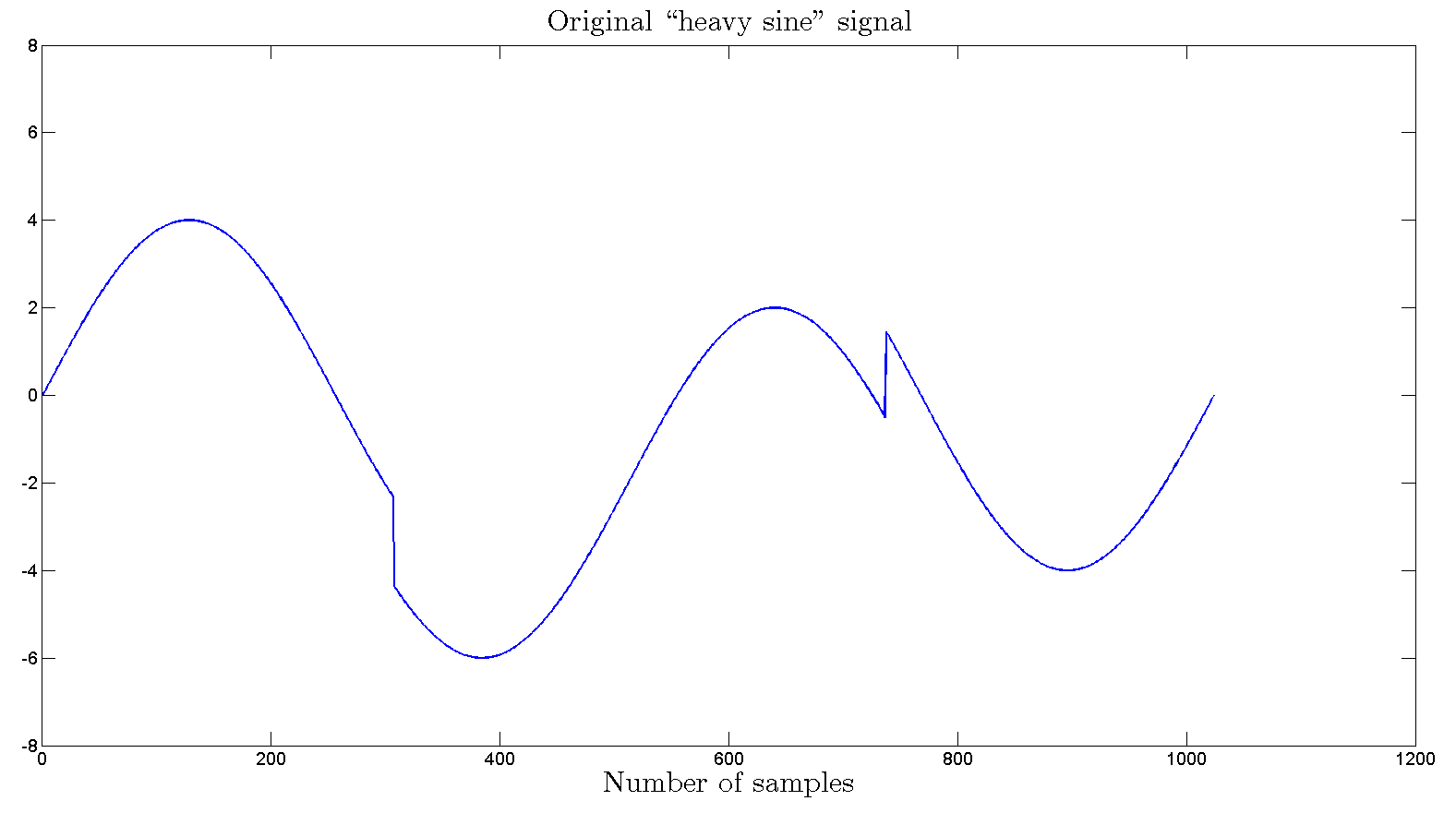}}
\subfigure[Noisy signal]
{\label{noisy}\includegraphics[width=75mm]{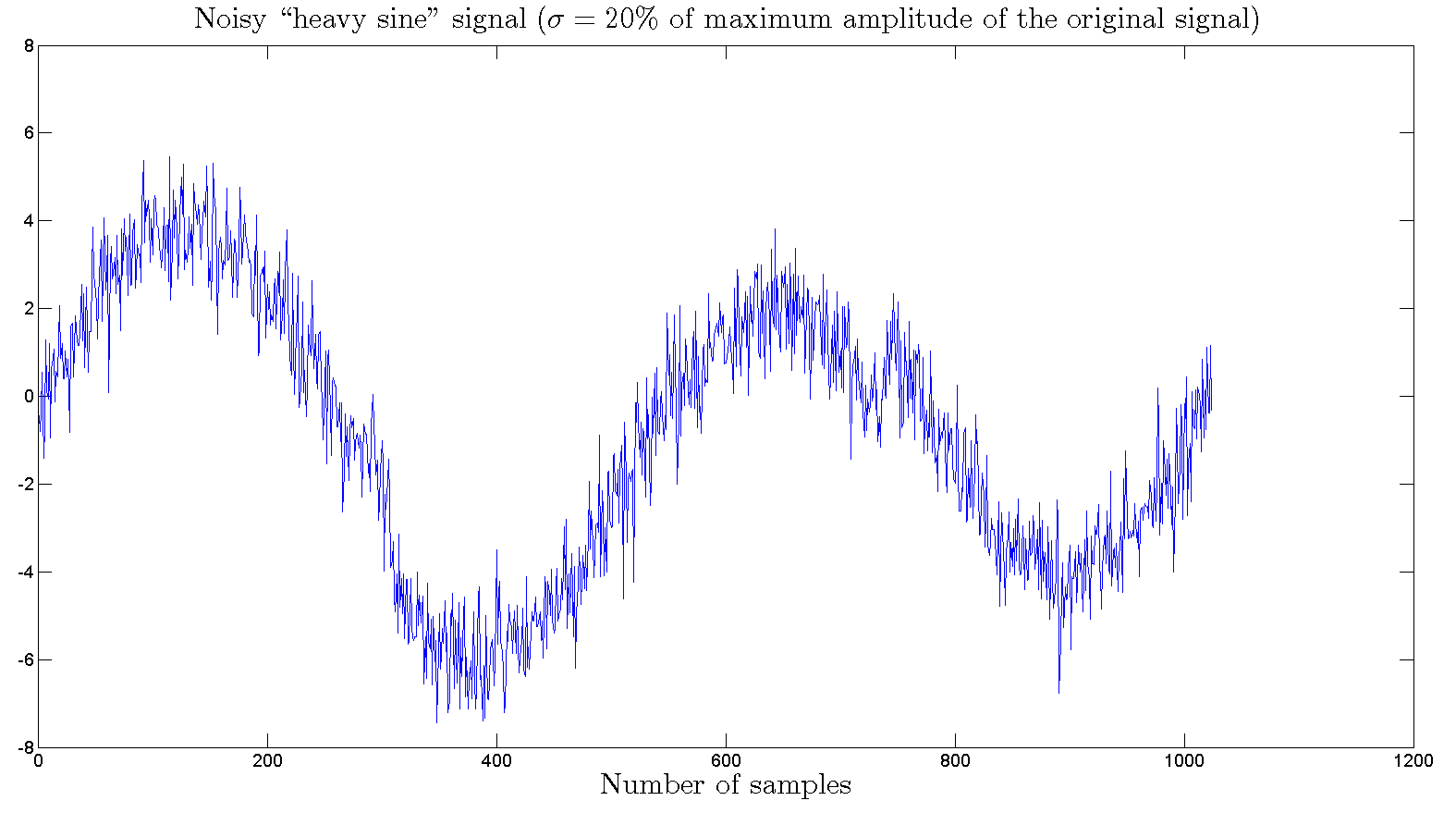}}
\subfigure[PES-$\ell_1$ with pyramid]
{\label{PESCDenoised}\includegraphics[width=75mm]{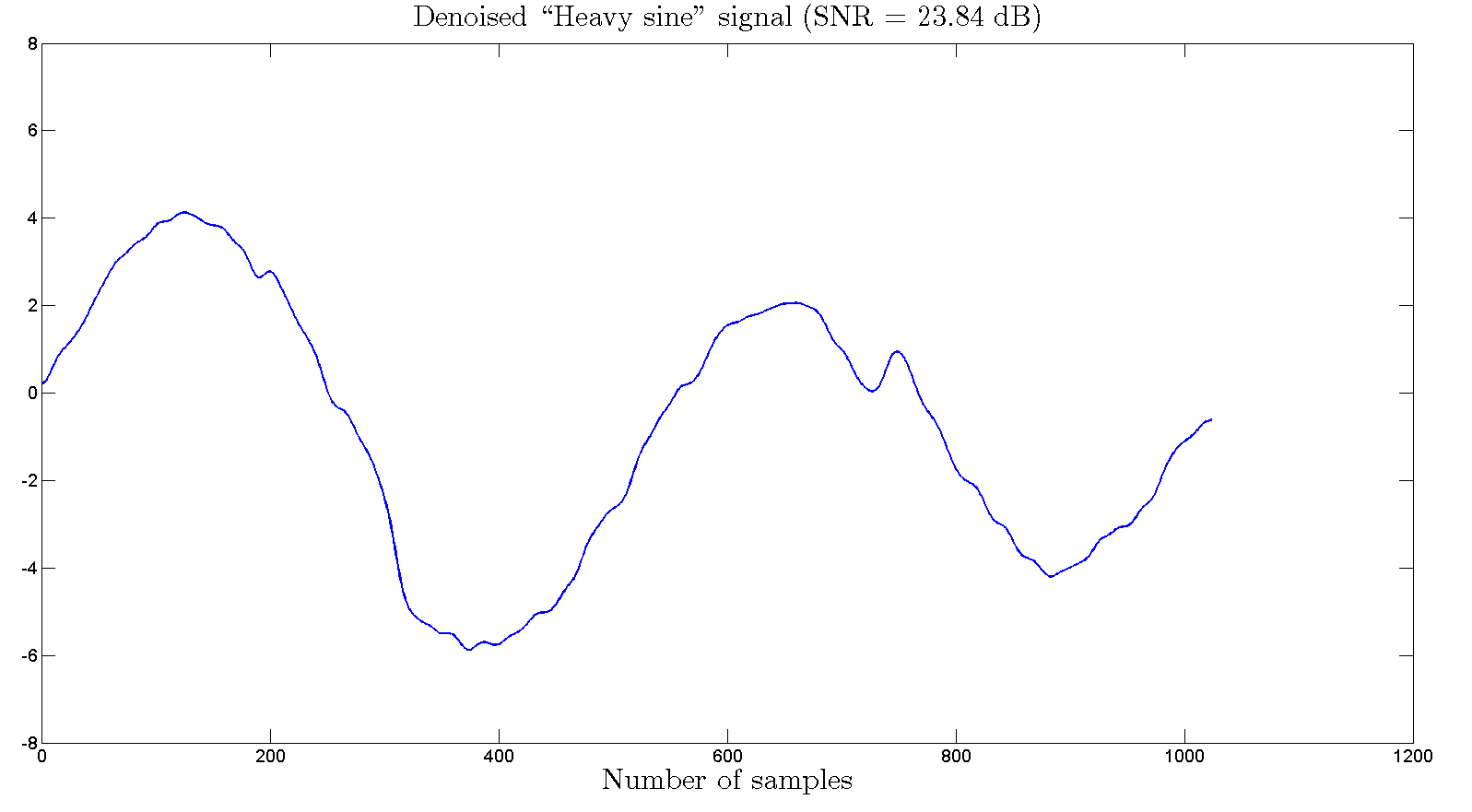}}
\subfigure[PES-$\ell_1$ with wavelet]
{\label{PESCDenoisedWavelet}\includegraphics[width=75mm]{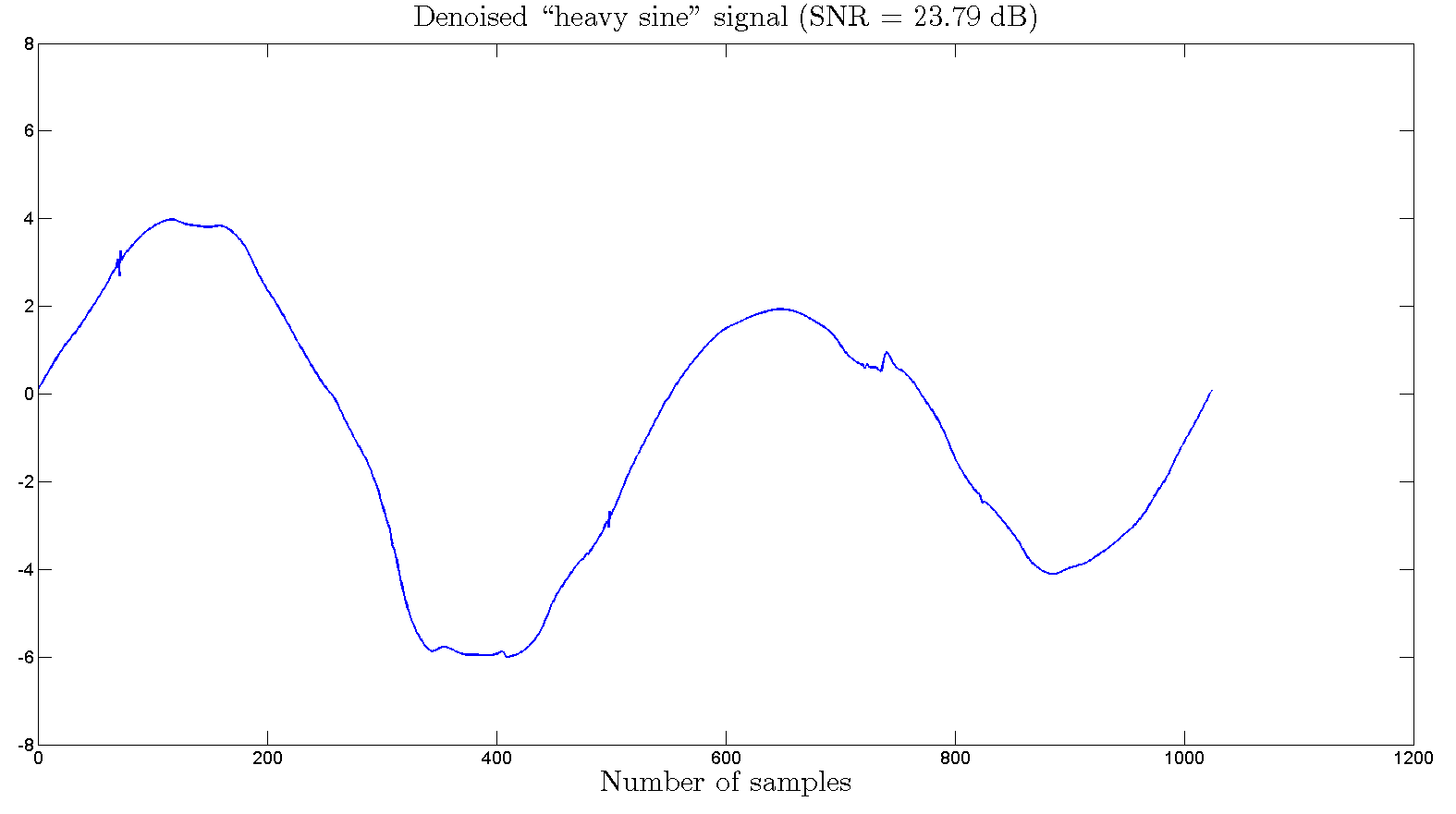}}
\subfigure[Wavelet denoising in MATLAB \cite{Rousseeuw,Mina}]
{\label{DenoisedWavelet}\includegraphics[width=75mm]{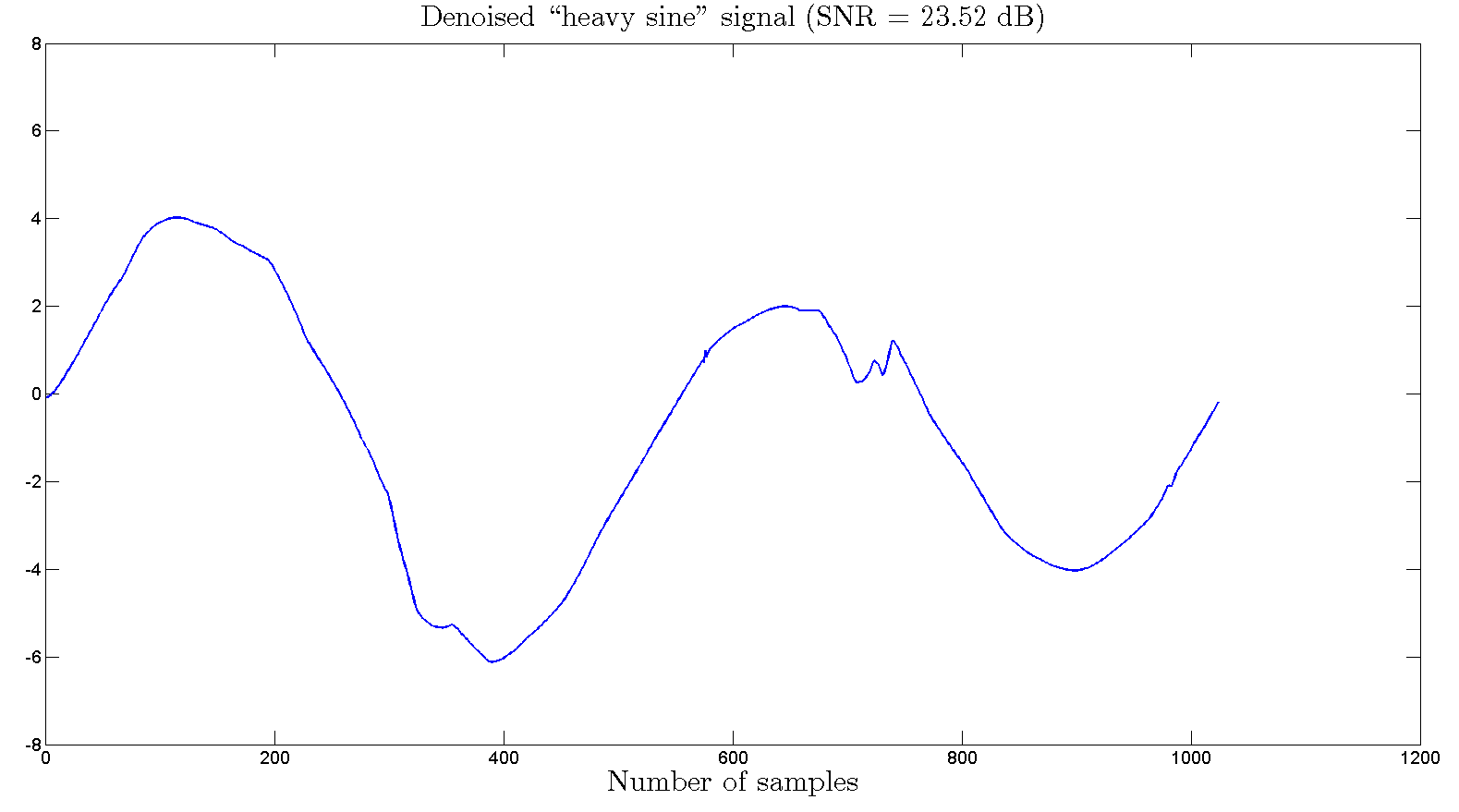}}
\subfigure[Wavelet denoising ``minimaxi" algorithm \cite{Donoho}]
{\label{Minimaxi}\includegraphics[width=75mm]{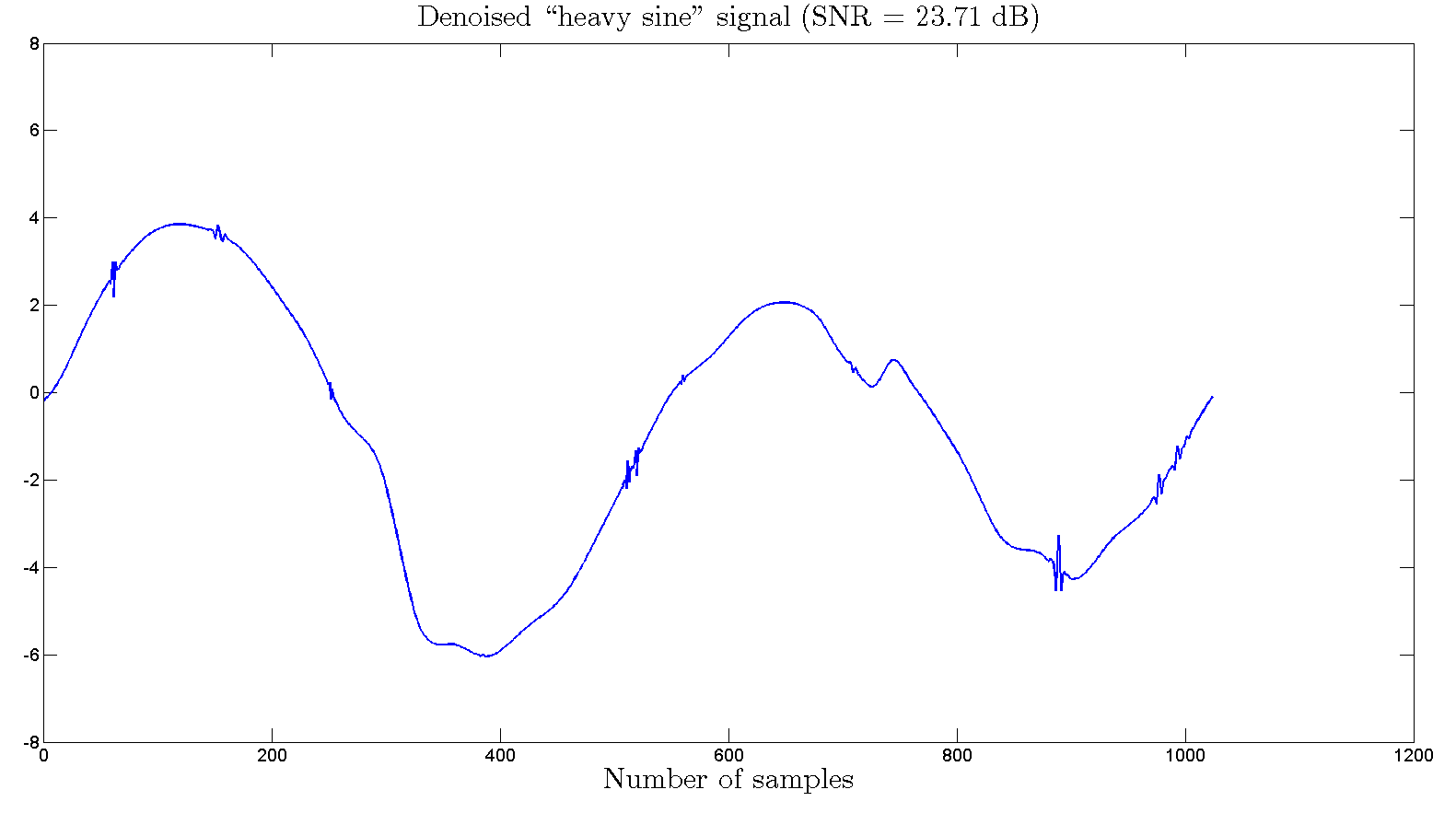}}
\subfigure[Wavelet denoising ``rigrsure" algorithm \cite{Johnstone}]
{\label{Rigrsure}\includegraphics[width=75mm]{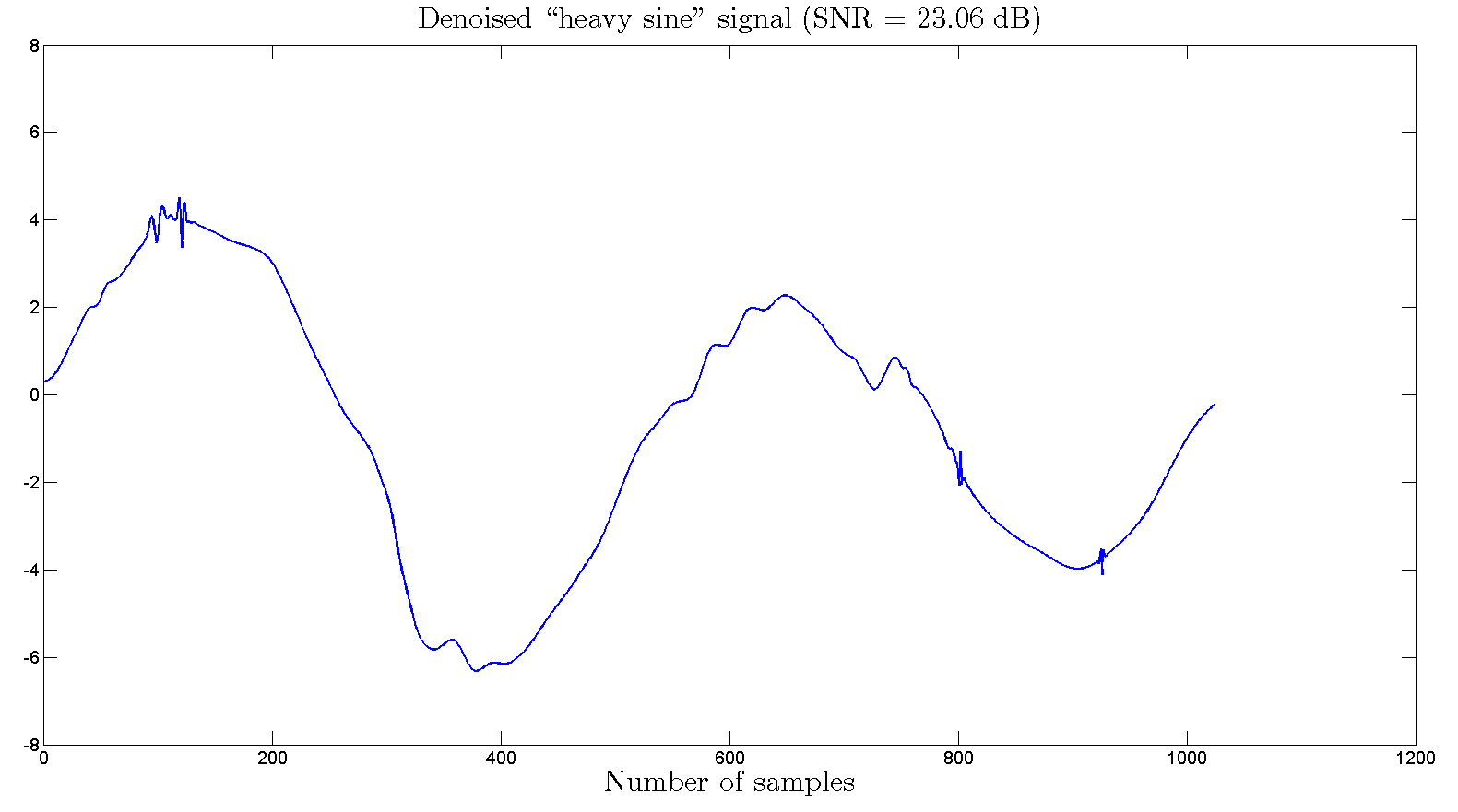}}
\subfigure[Wavelet denoising with T = 3$\hat{\sigma}$ \cite{Donoho,FowlerDen}]
{\label{Soft-thresholding}\includegraphics[width=75mm]{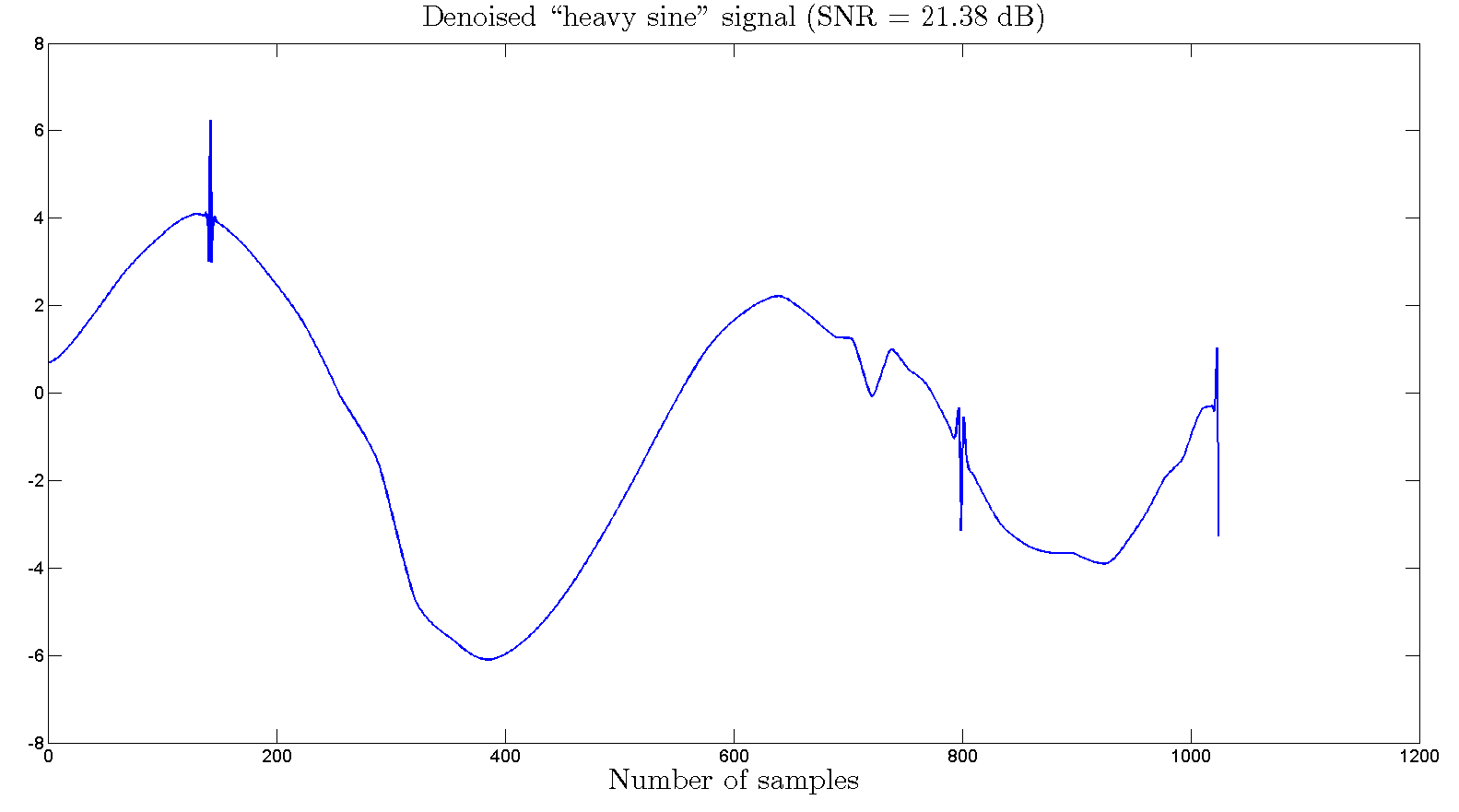}}
\caption{(a) Original ``heavy sine" signal, (b) signal corrupted with Gaussian noise with $\sigma = 20 \%$ of maximum amplitude of the original signal, and denoised signal using (c) PES-$\ell_1$-ball with pyramid; SNR = 23.84 dB and, (d) PES-$\ell_1$-ball with wavelet; SNR = 23.79 dB, (e) Wavelet denoising in Matlab; SNR = 23.52 dB \cite{Rousseeuw,Mina}, (f) Wavelet denoising ``minimaxi" algorithm \cite{Donoho}; SNR = 23.71 dB, (g) Wavelet denoising ``rigrsure" algorithm \cite{Johnstone}; SNR = 23.06 dB, (h) Wavelet denoising with T = 3$\hat{\sigma}$ \cite{Donoho,FowlerDen}; SNR = 21.38 dB.}
\label{app:CUSPDenoising}
\end{figure*}

\begin{figure*}[htb!]
\centering
\subfigure[Original signal]
{\label{original}\includegraphics[width=75mm]{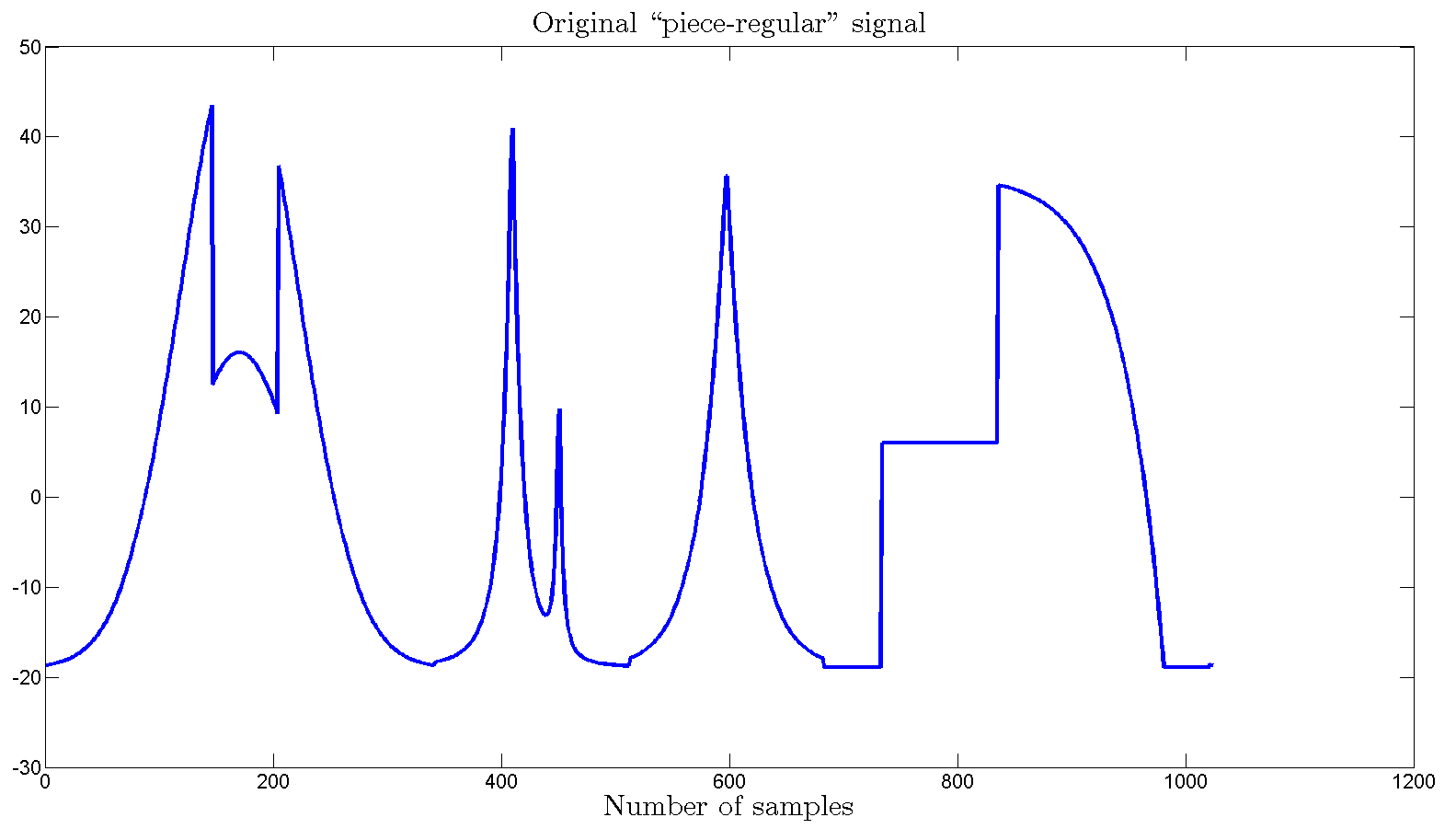}}
\subfigure[Noisy signal]
{\label{noisy}\includegraphics[width=75mm]{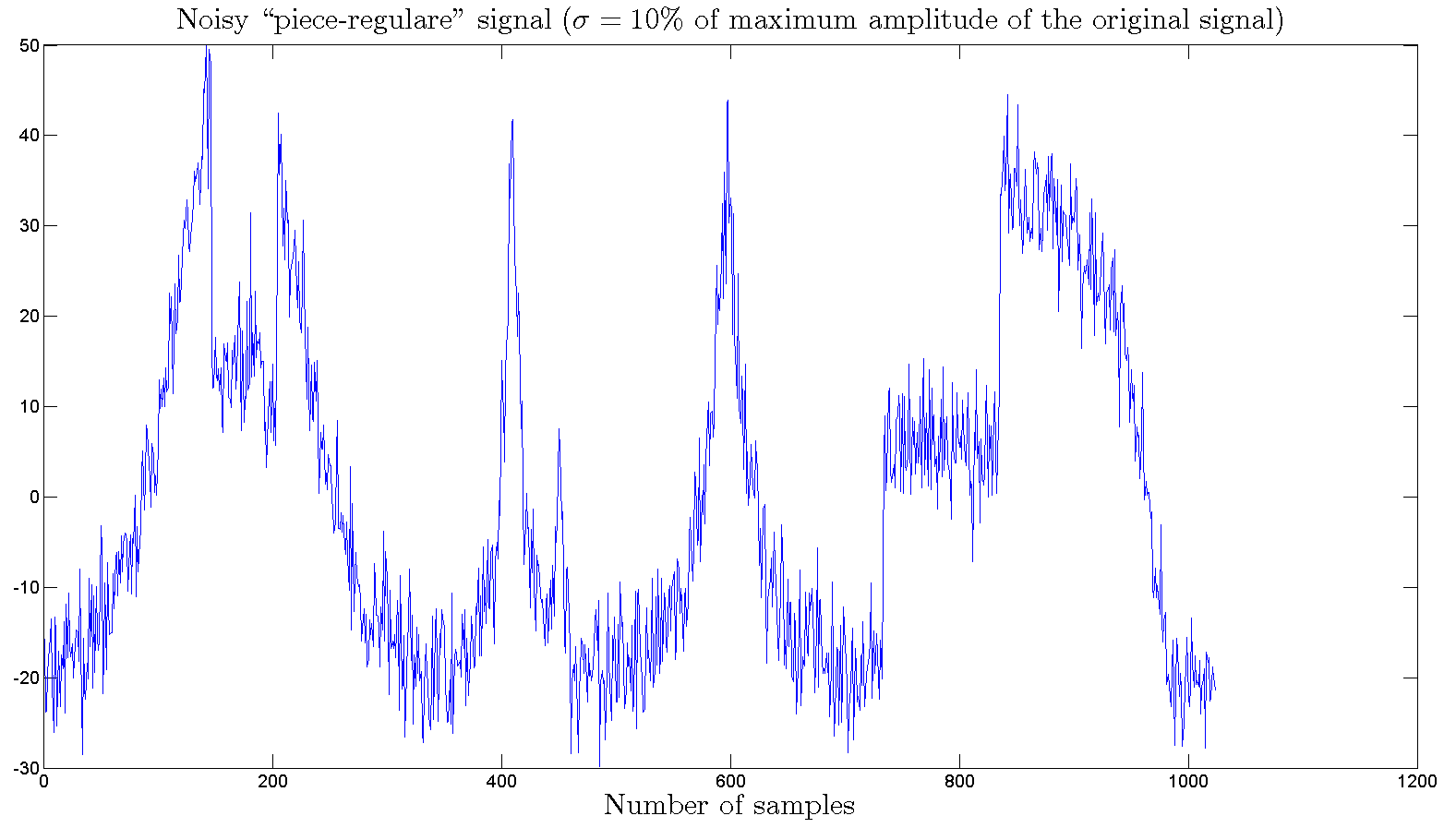}}
\subfigure[PES-$\ell_1$ with pyramid]
{\label{PESCDenoised}\includegraphics[width=75mm]{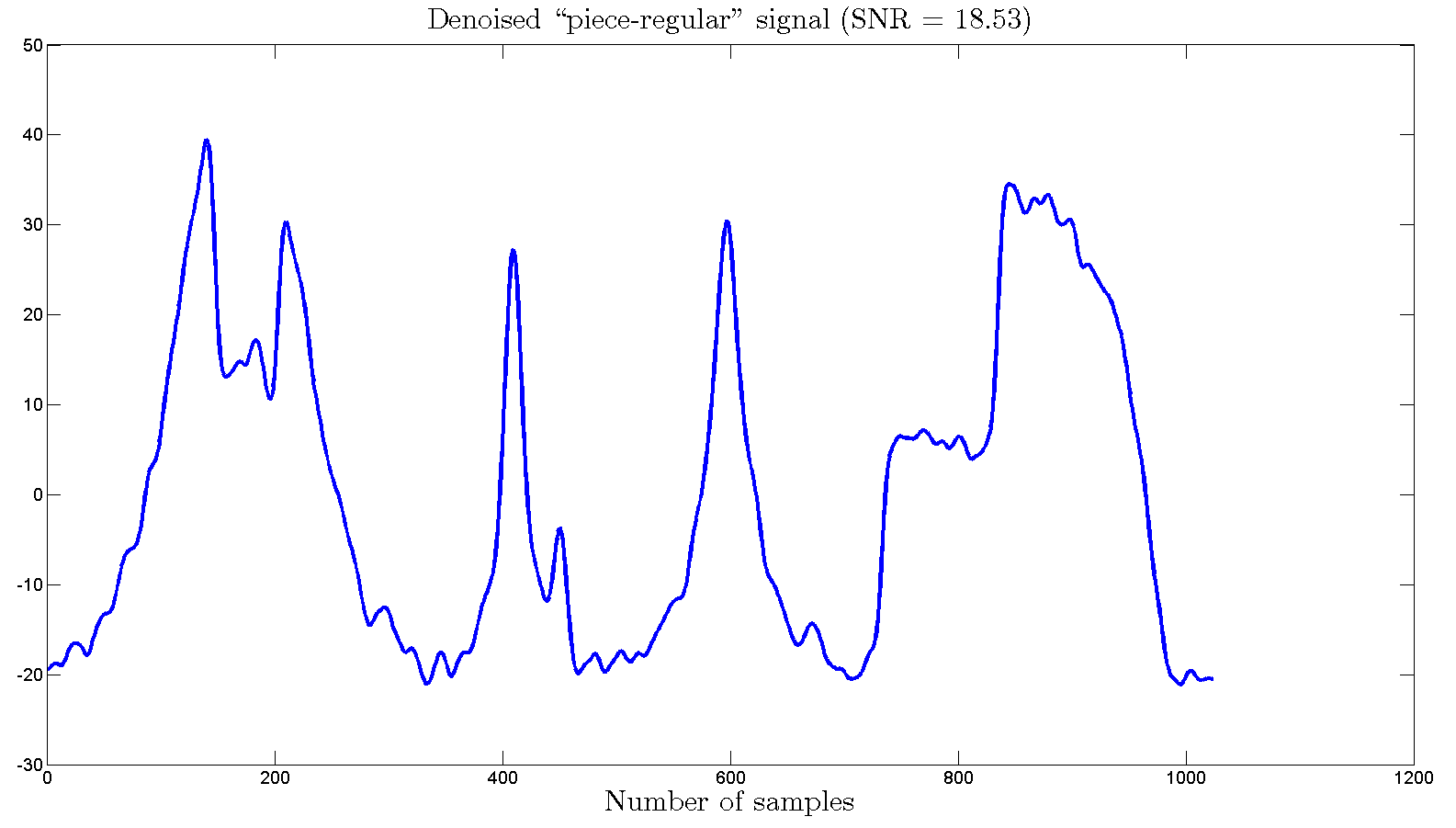}}
\subfigure[PES-$\ell_1$ with wavelet]
{\label{PESCDenoisedWavelet}\includegraphics[width=75mm]{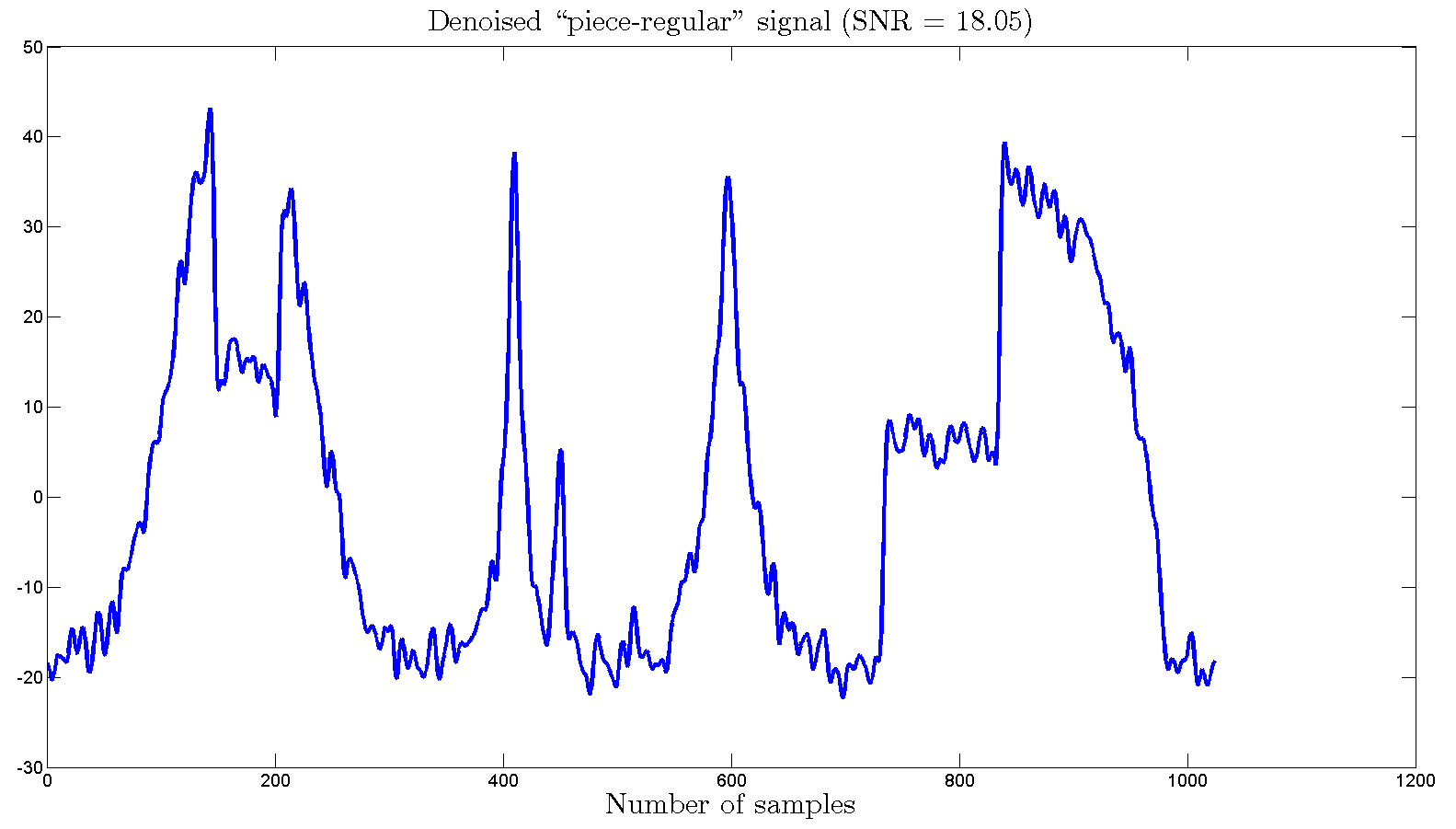}}
\subfigure[Wavelet denoising in MATLAB \cite{Rousseeuw,Mina}]
{\label{DenoisedWavelet}\includegraphics[width=75mm]{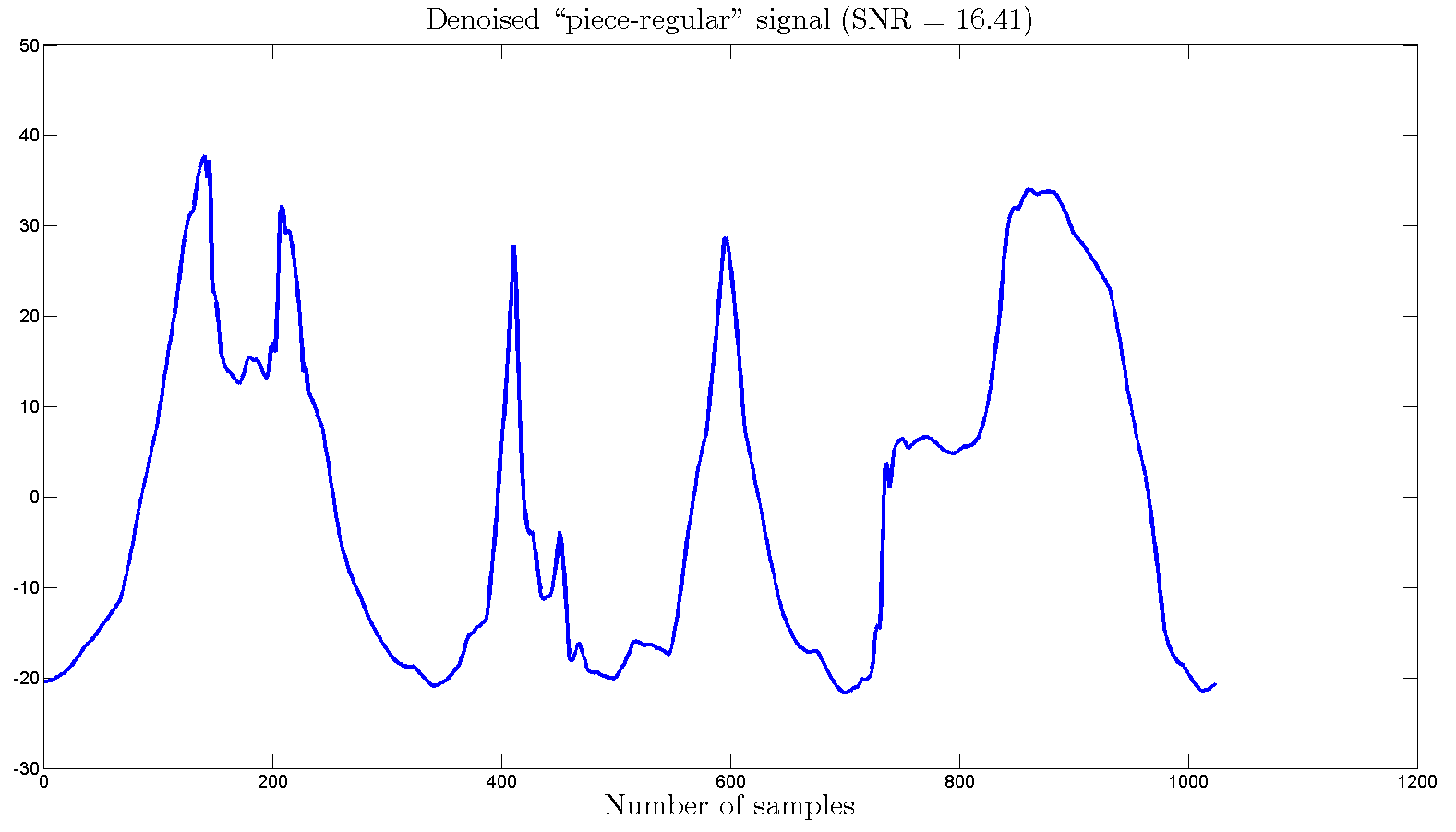}}
\subfigure[Wavelet denoising ``minimaxi" algorithm \cite{Donoho}]
{\label{Minimaxi}\includegraphics[width=75mm]{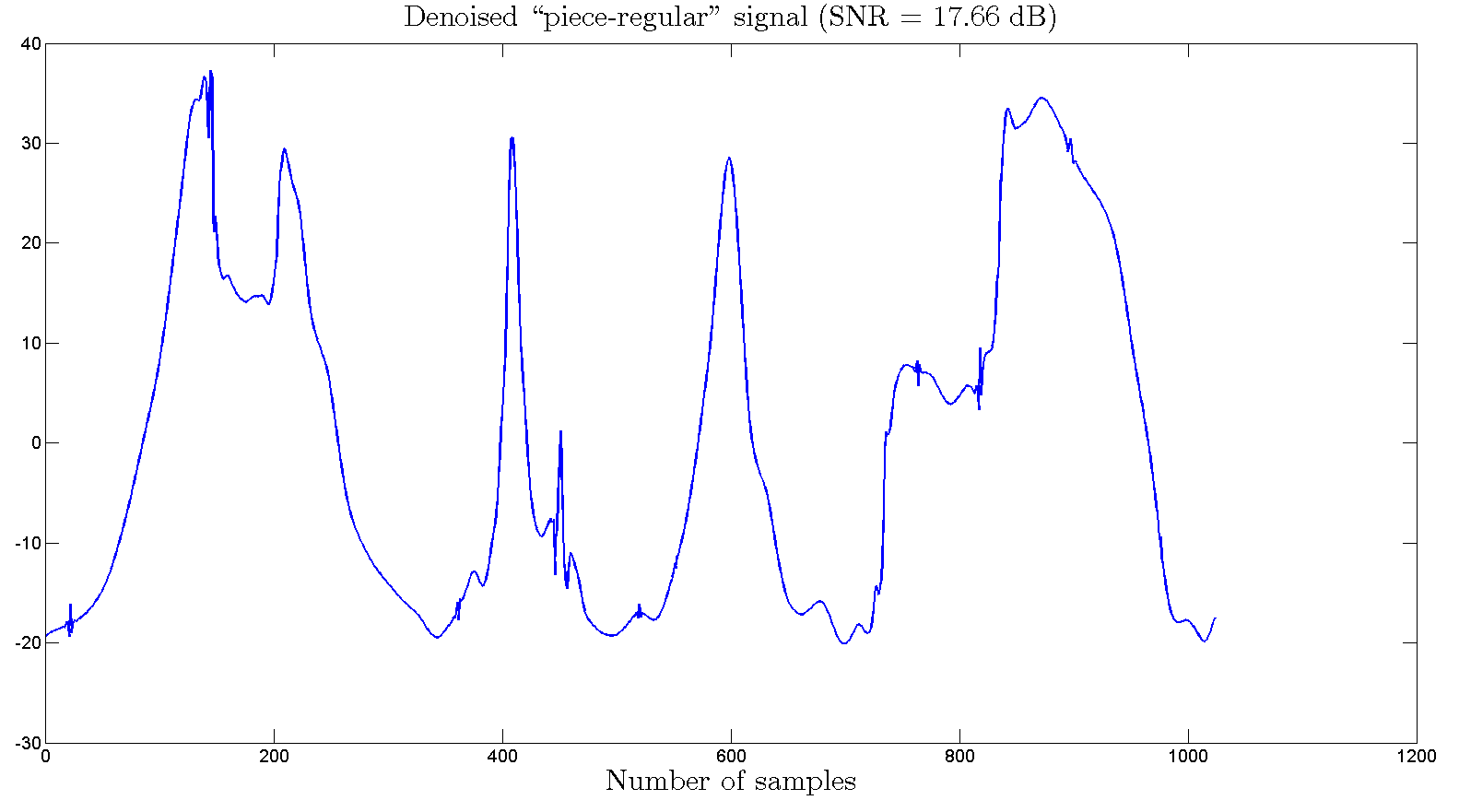}}
\subfigure[Wavelet denoising ``rigrsure" algorithm \cite{Johnstone}]
{\label{Rigrsure}\includegraphics[width=75mm]{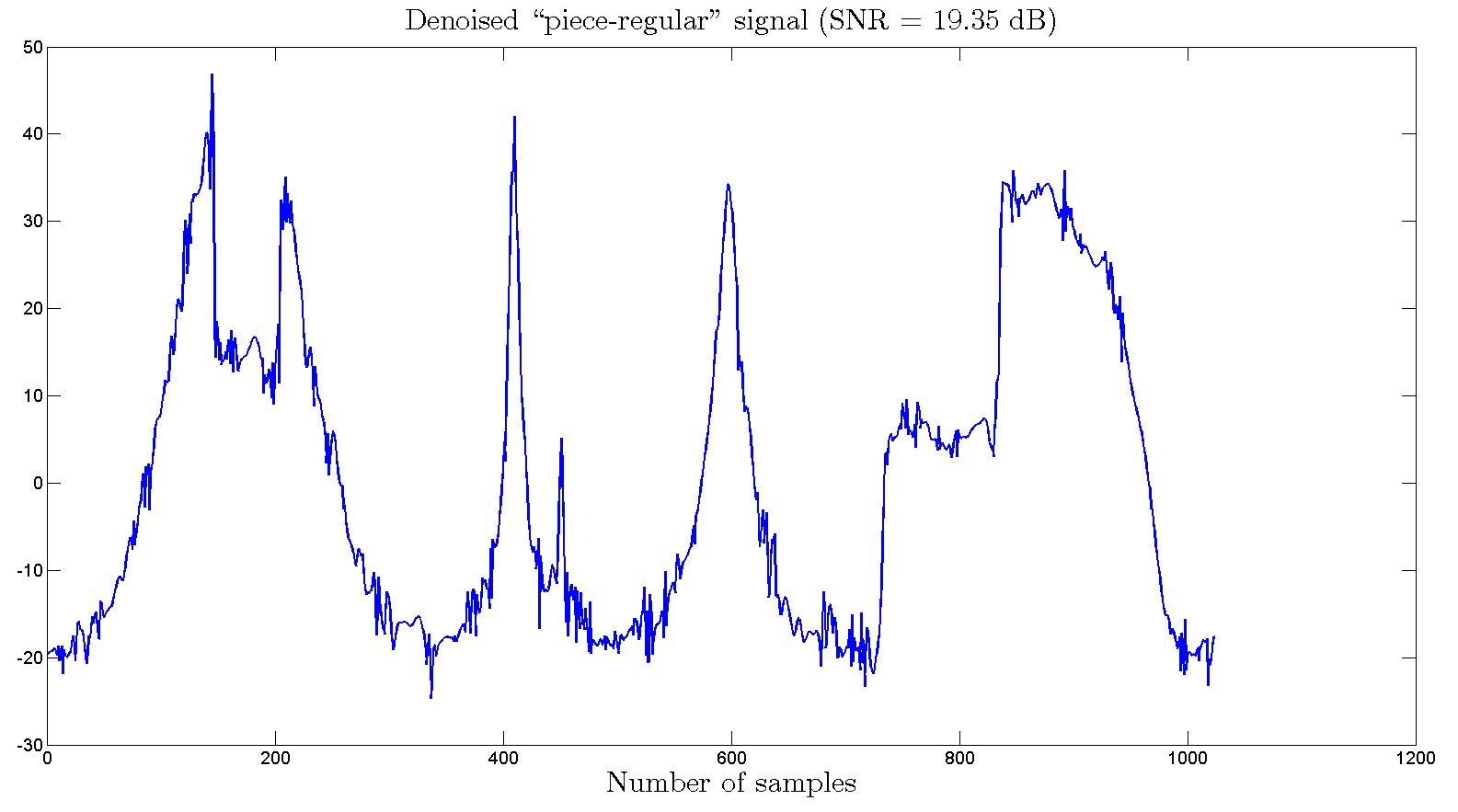}}
\subfigure[Wavelet denoising with T = 3$\hat{\sigma}$ \cite{Donoho,FowlerDen}]
{\label{Soft-thresholding}\includegraphics[width=75mm]{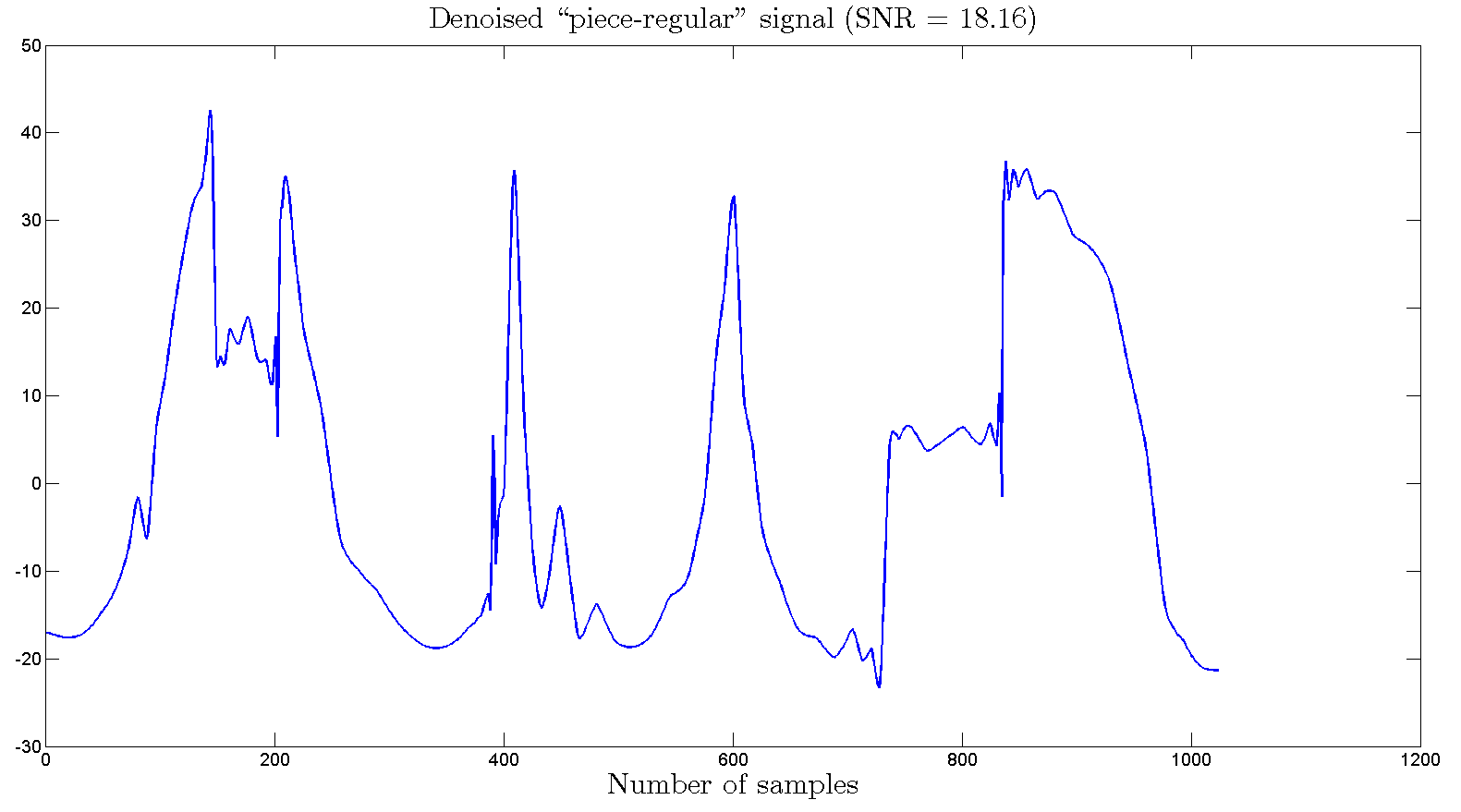}}
\caption{(a) Original ``cusp" signal, (b) signal corrupted with Gaussian noise with $\sigma = 10 \%$ of maximum amplitude of the original signal, and denoised signal using (c) PES-$\ell_1$-ball with pyramid; SNR = 23.84 dB and, (d) PES-$\ell_1$-ball with wavelet; SNR = 23.79 dB, (e) Wavelet denoising in Matlab; SNR = 23.52 dB \cite{Rousseeuw,Mina}, (f) Wavelet denoising ``minimaxi" algorithm \cite{Donoho}; SNR = 23.71 dB, (g) Wavelet denoising ``rigrsure" algorithm \cite{Johnstone}; SNR = 23.06 dB, (h) Wavelet denoising with T = 3$\hat{\sigma}$ \cite{Donoho,FowlerDen}; SNR = 21.38 dB.}
\label{app:CUSPDenoising}
\end{figure*}

\begin{figure*}[htb!]
\centering
\subfigure[Original signal]
{\label{original}\includegraphics[width=75mm]{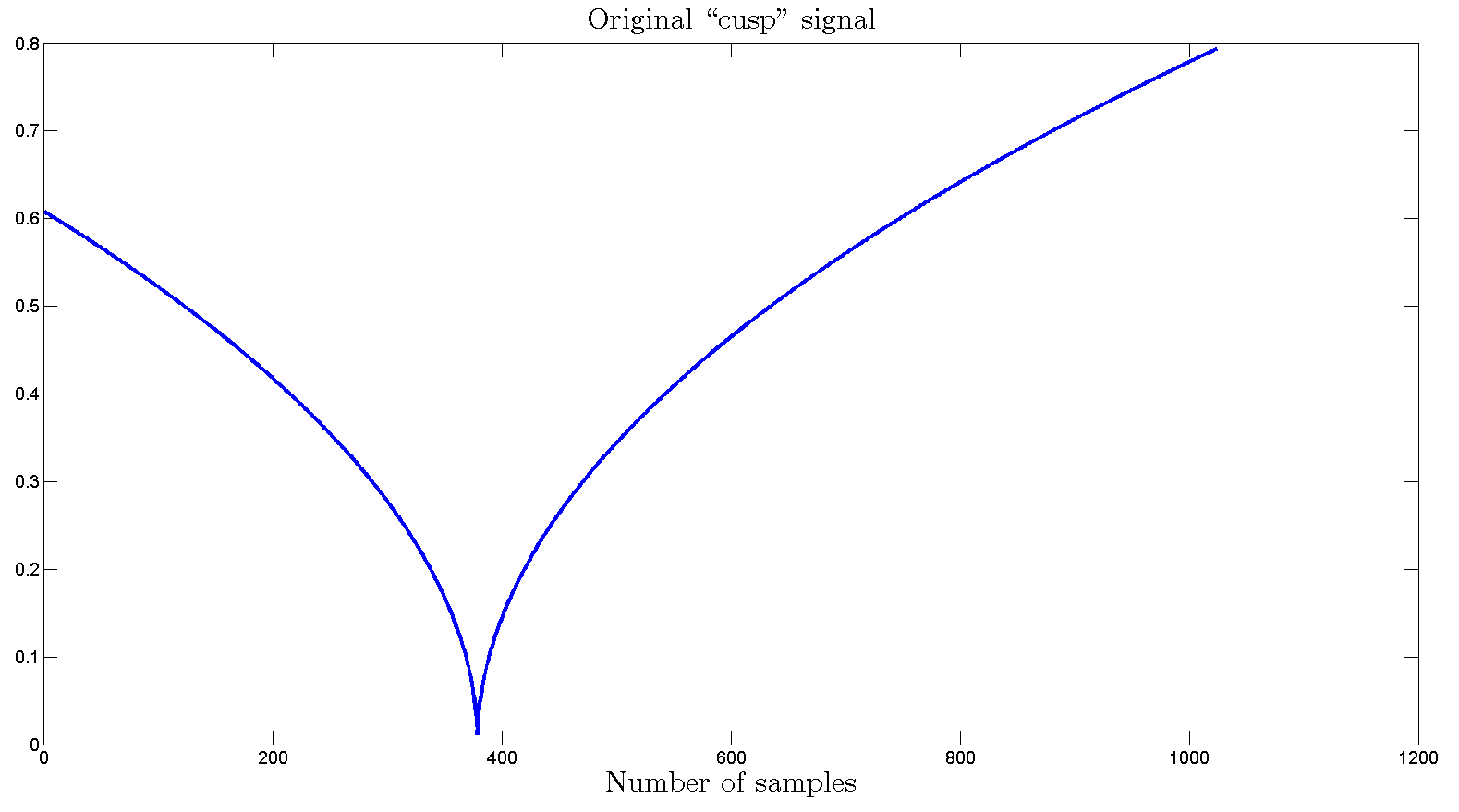}}
\subfigure[Noisy signal]
{\label{noisy}\includegraphics[width=75mm]{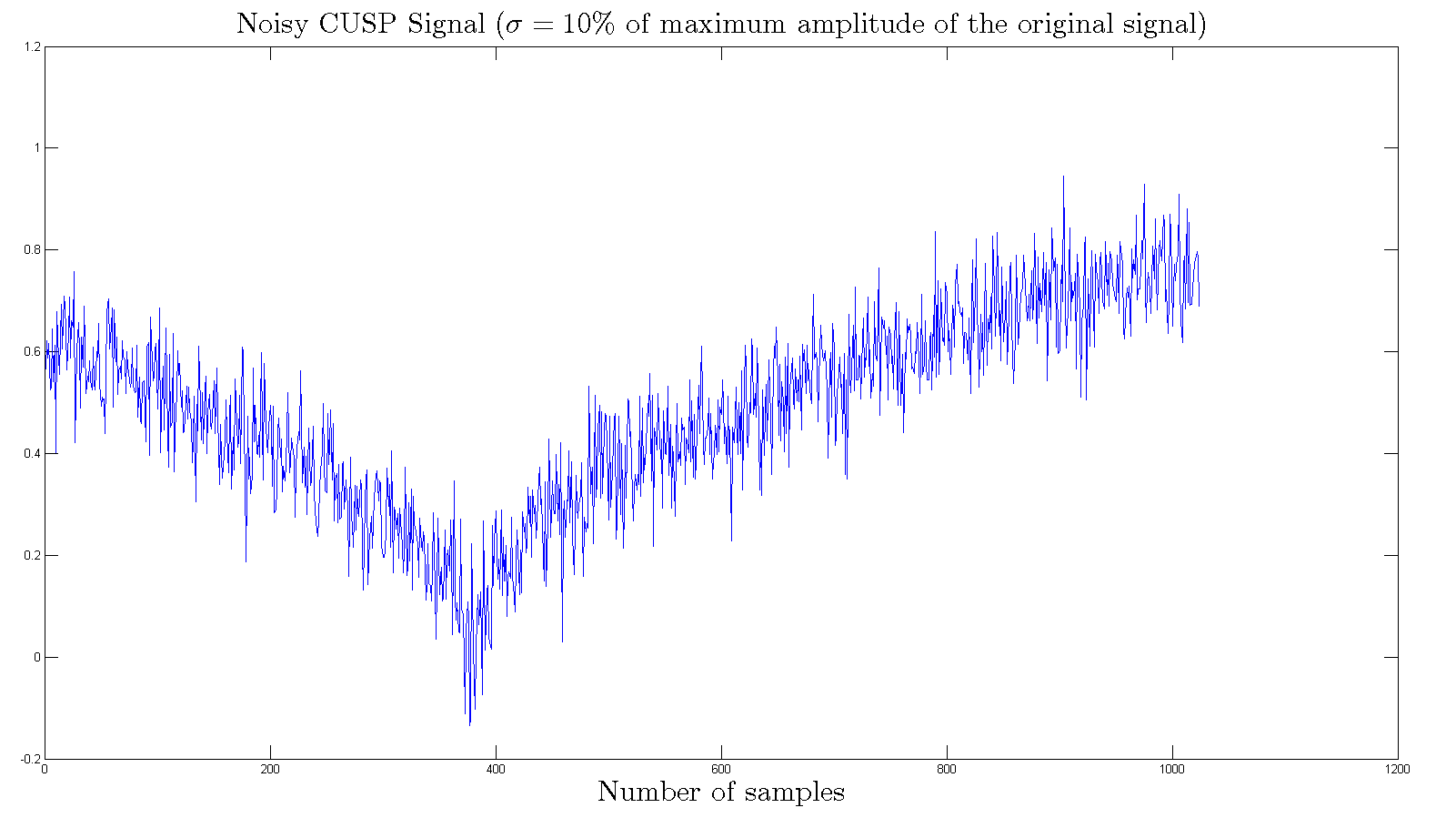}}
\subfigure[PES-$\ell_1$ with pyramid]
{\label{PESCDenoised}\includegraphics[width=75mm]{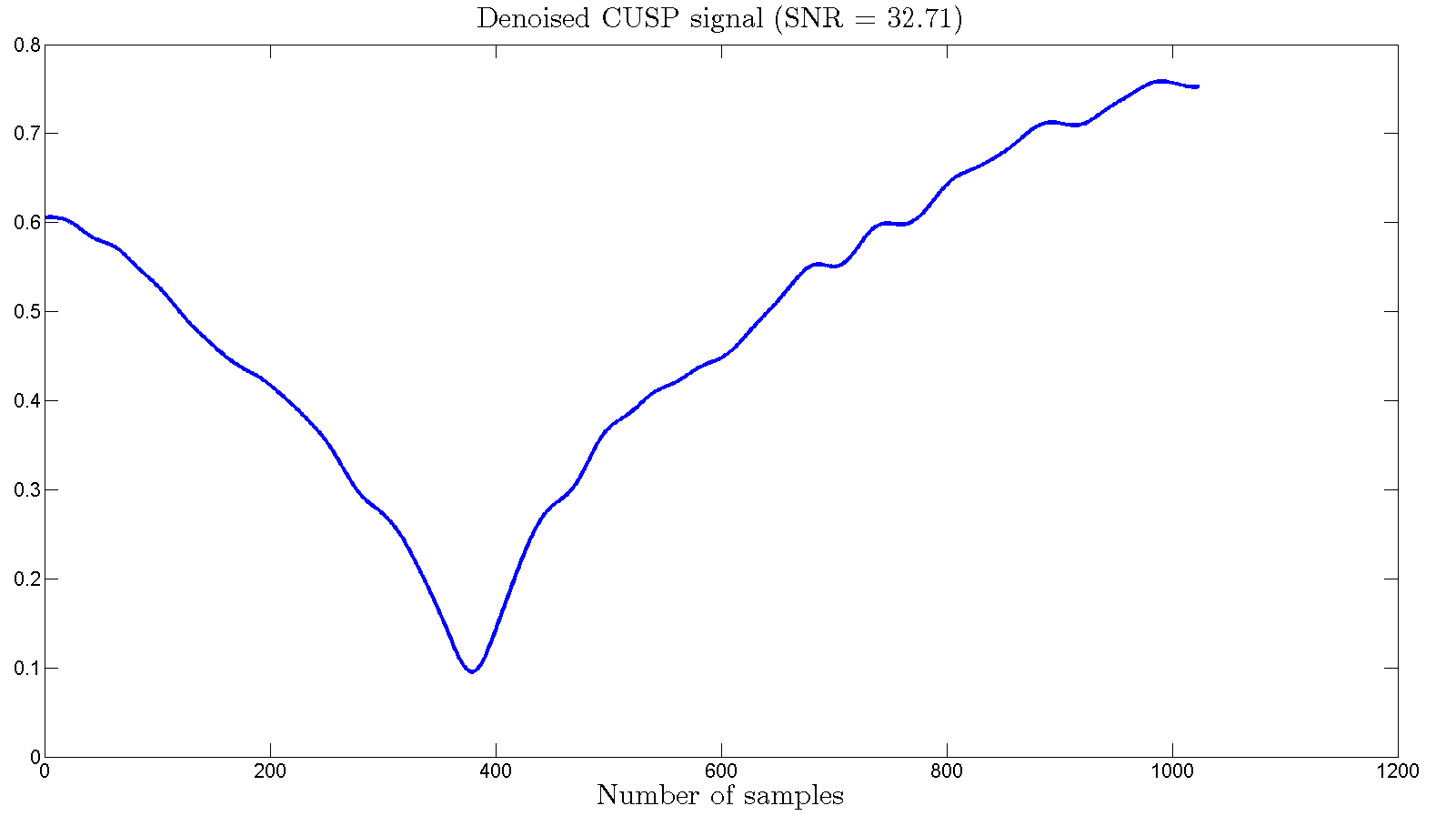}}
\subfigure[PES-$\ell_1$ with wavelet]
{\label{PESCDenoisedWavelet}\includegraphics[width=75mm]{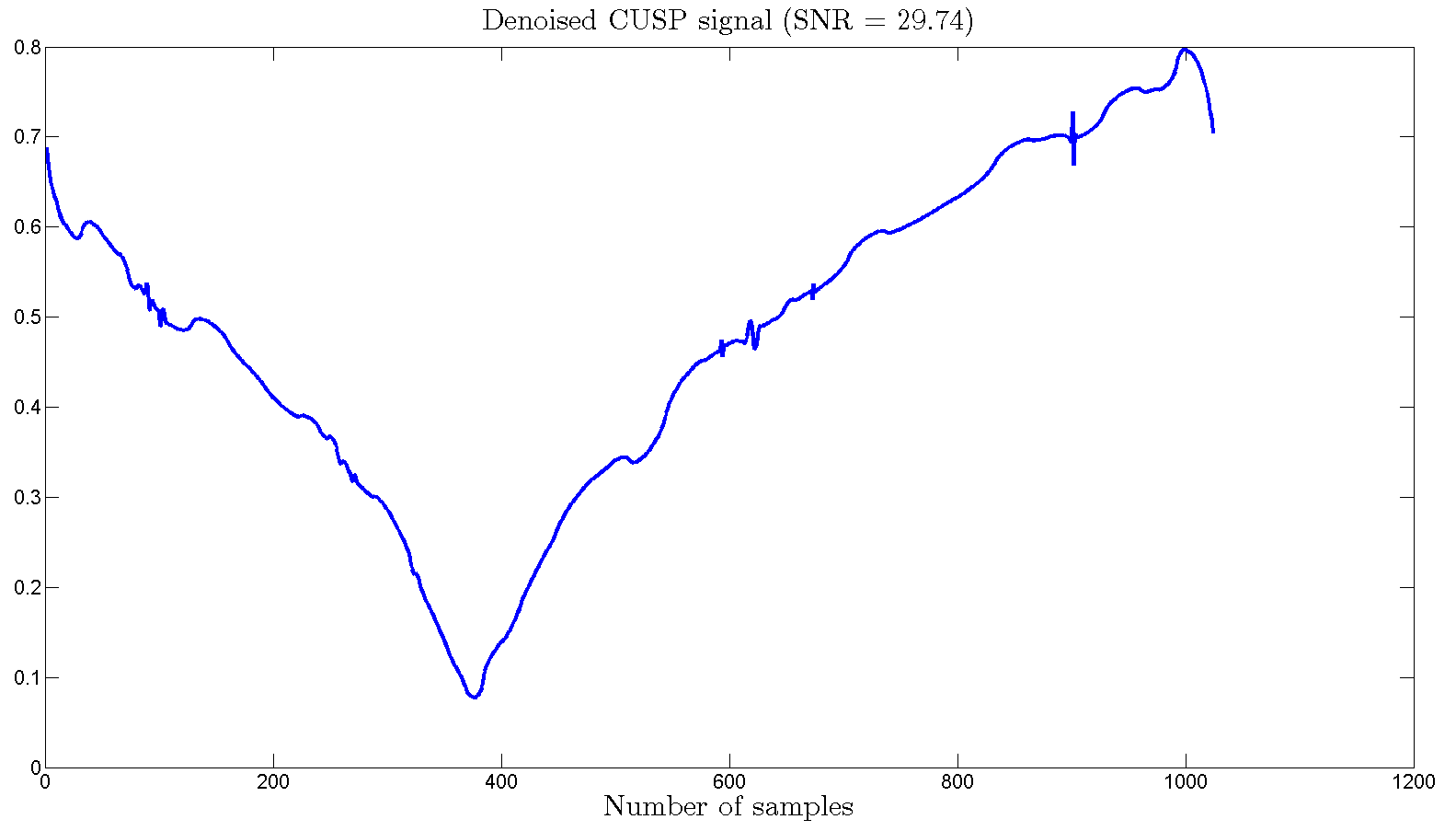}}
\subfigure[Wavelet denoising in MATLAB \cite{Rousseeuw,Mina}]
{\label{DenoisedWavelet}\includegraphics[width=75mm]{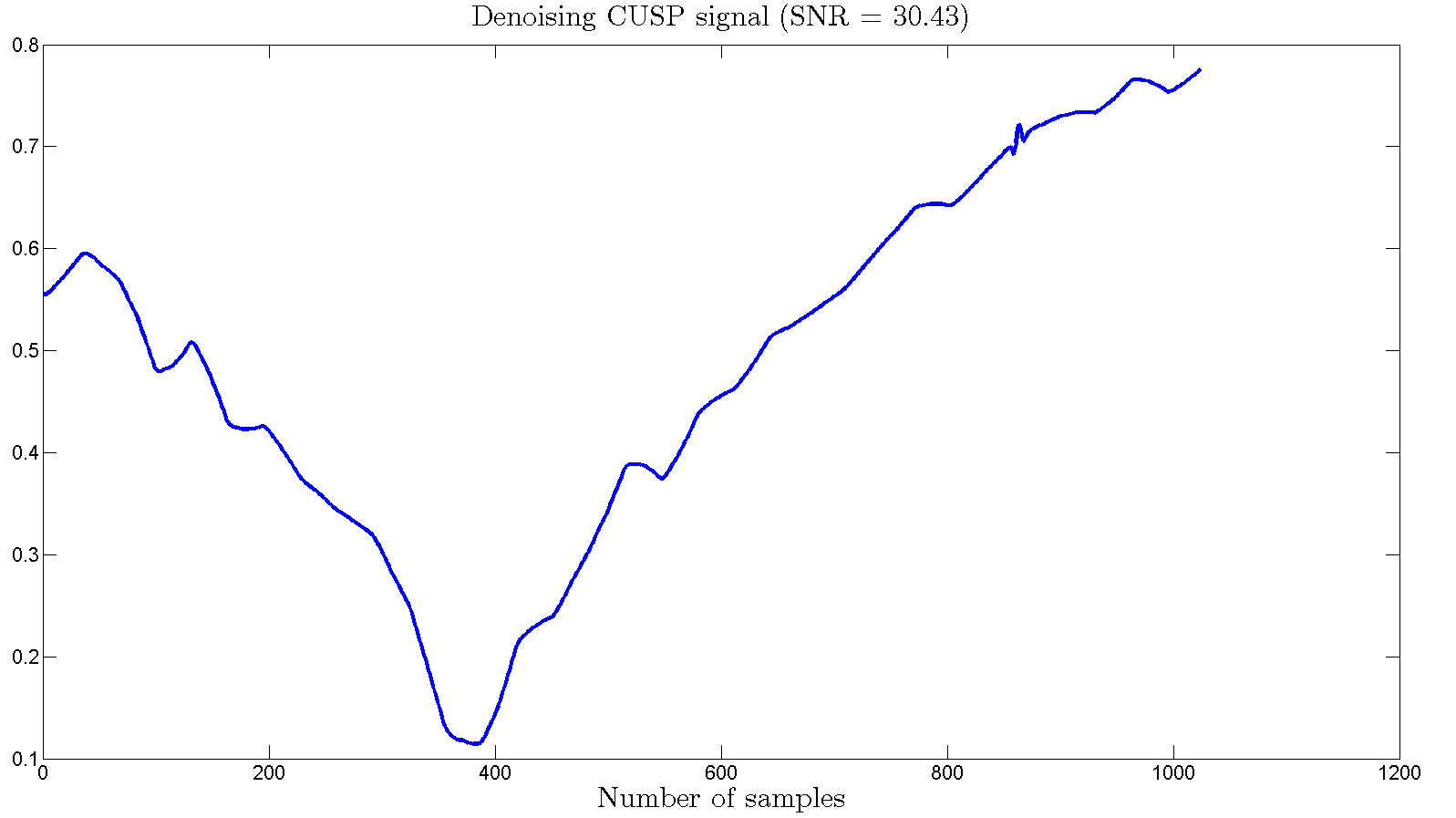}}
\subfigure[Wavelet denoising ``minimaxi" algorithm \cite{Donoho}]
{\label{Minimaxi}\includegraphics[width=75mm]{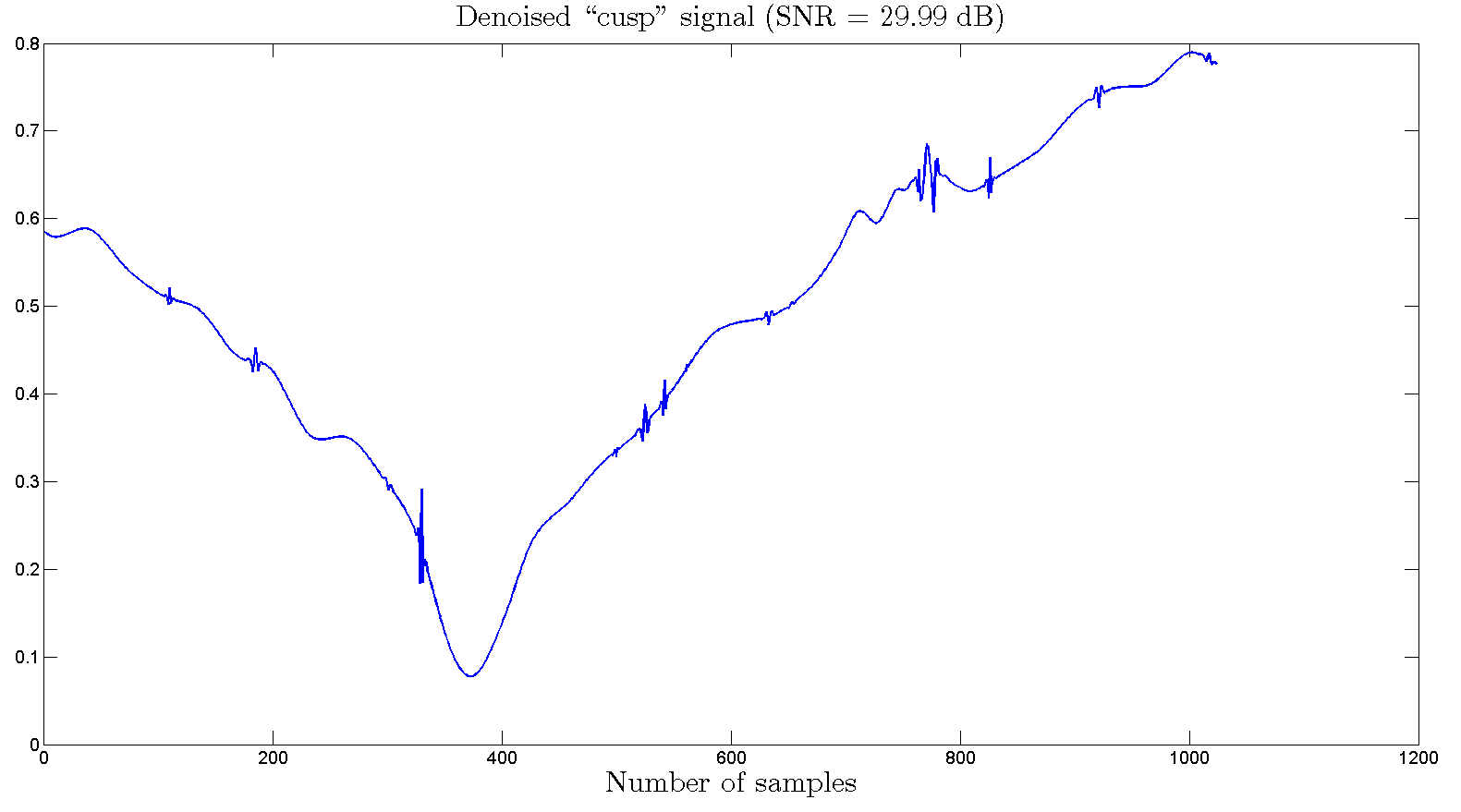}}
\subfigure[Wavelet denoising ``rigrsure" algorithm \cite{Johnstone}]
{\label{Rigrsure}\includegraphics[width=75mm]{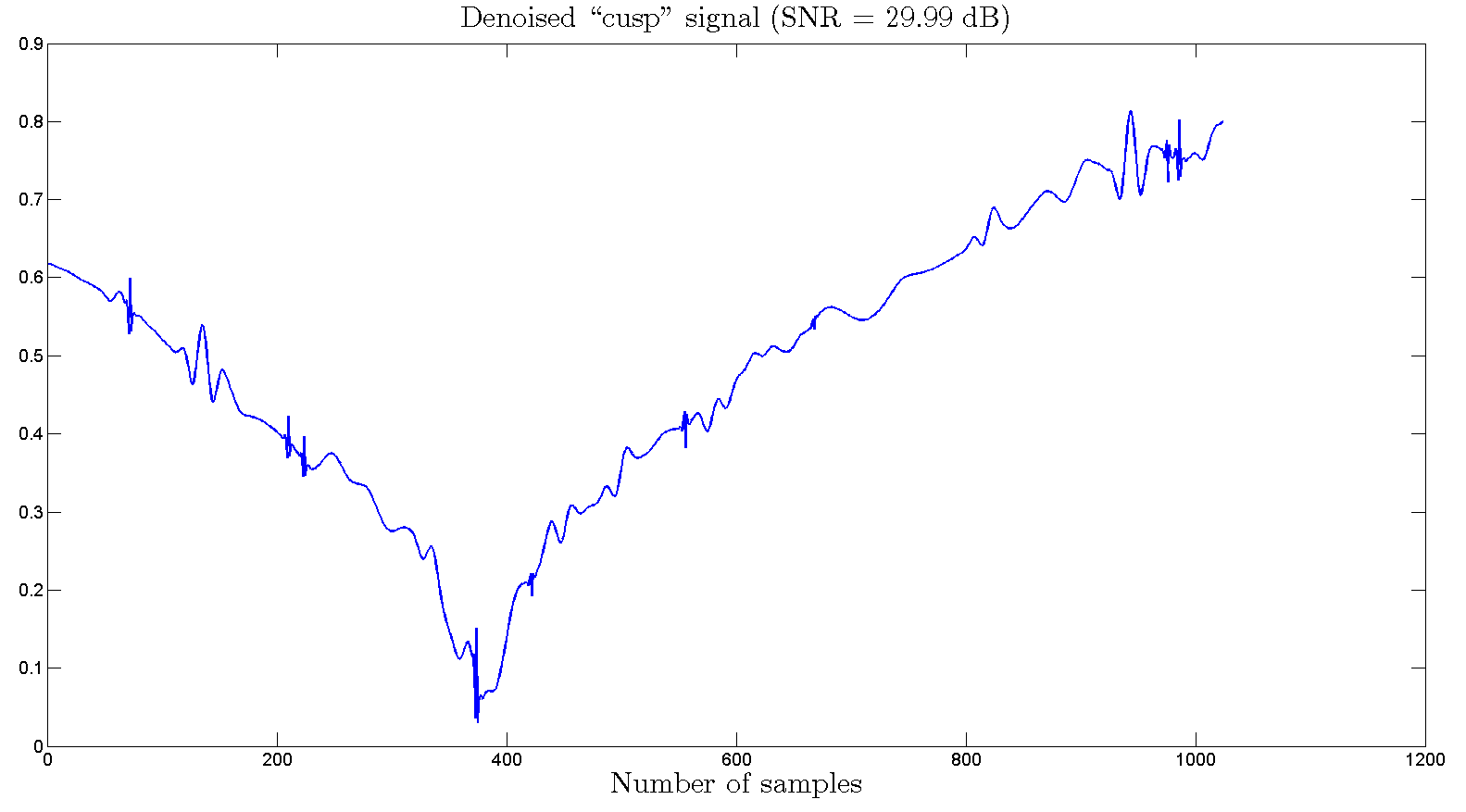}}
\subfigure[Wavelet denoising with T = 3$\hat{\sigma}$ \cite{Donoho,FowlerDen}]
{\label{Soft-thresholding}\includegraphics[width=75mm]{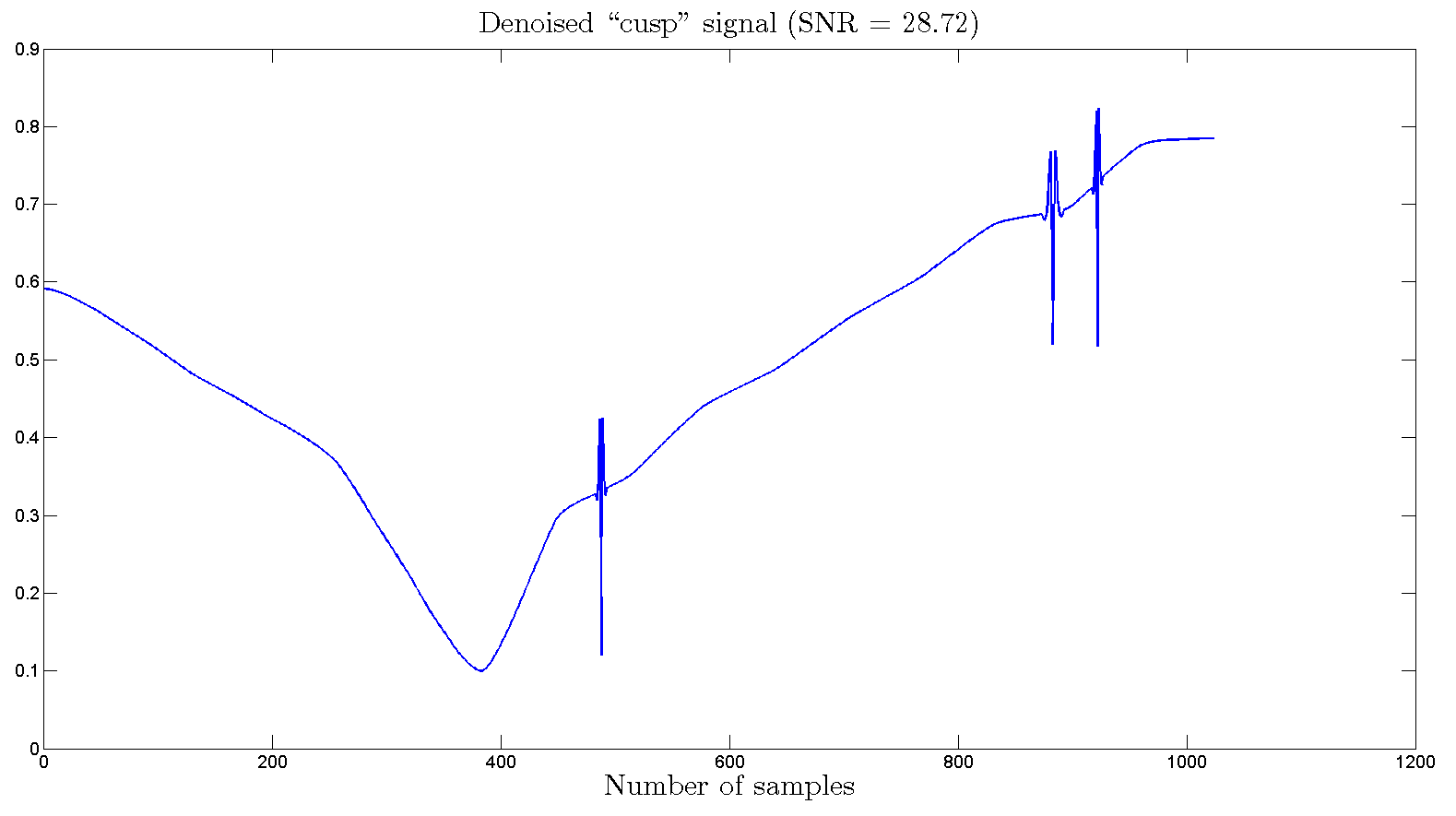}}
\caption{(a) Original ``cusp" signal, (b) signal corrupted with Gaussian noise with $\sigma = 10 \%$ of maximum amplitude of the original signal, and denoised signal using (c) PES-$\ell_1$-ball with pyramid; SNR = 23.84 dB and, (d) PES-$\ell_1$-ball with wavelet; SNR = 23.79 dB, (e) Wavelet denoising in Matlab; SNR = 23.52 dB \cite{Rousseeuw,Mina}, (f) Wavelet denoising ``minimaxi" algorithm \cite{Donoho}; SNR = 23.71 dB, (g) Wavelet denoising ``rigrsure" algorithm \cite{Johnstone}; SNR = 23.06 dB, (h) Wavelet denoising with T = 3$\hat{\sigma}$ \cite{Donoho,FowlerDen}; SNR = 21.38 dB.}
\label{app:CUSPDenoising}
\end{figure*}

\section{SIMULATION RESULTS}
\label{sec:SIMULATION RESULTS}
Epigraph set based threshold selection is compared with wavelet denoising methods used in MATLAB \cite{Rousseeuw,Mina,vetterli,Donoho}. The ``heavy sine" signal shown in Fig. \ref{original} is corrupted by a zero mean Gaussian noise with $\sigma = 20 \%$ of the maximum amplitude of the original signal. The signal is restored using PES-$\ell_1$ with pyramid structure, PES-$\ell_1$ with wavelet, MATLAB's wavelet multivariate denoising algorithm \cite{Rousseeuw,Mina}, MATLAB's soft-thresholding denoising algorithm (for ``minimaxi" and ``rigrsure" thresholds), and wavelet thresholding denoising method. The denoised signals are shown in Fig.~\ref{PESCDenoised}, \ref{PESCDenoisedWavelet}, \ref{DenoisedWavelet}, \ref{Minimaxi}, \ref{Rigrsure}, and \ref{Soft-thresholding} with SNR values equal to 23.84, 23.79, 23.52, 23.71, 23.06 dB, and 21.38, respectively. On the average, the proposed PES-$\ell_1$ with pyramid and  PES-$\ell_1$ with wavelet method produce better thresholds than the other soft-thresholding methods. MATLAB codes of the denoising algorithms and other simulation examples are available in the following web-page: http://signal.ee.bilkent.edu.tr/1DDenoisingSoftware.html.

Results for other test signals in MATLAB are presented in Table \ref{tab:1}. These results are obtained by averaging the SNR values after repeating the simulations for 300 times. The SNR is calculated using the formula: $\textrm{SNR} = 20\times log_{10}(\|\textbf{w}_{orig}\| / \|\textbf{w}_{orig}-\textbf{w}_{rec}\|)$. In this lecture note, it is shown that soft-denoising threshold can be determined using basic linear algebra.

\begin{figure}[ht!]
\centering
\subfigure[Signal 1]
{\label{Signal 3}\includegraphics[width=75mm]{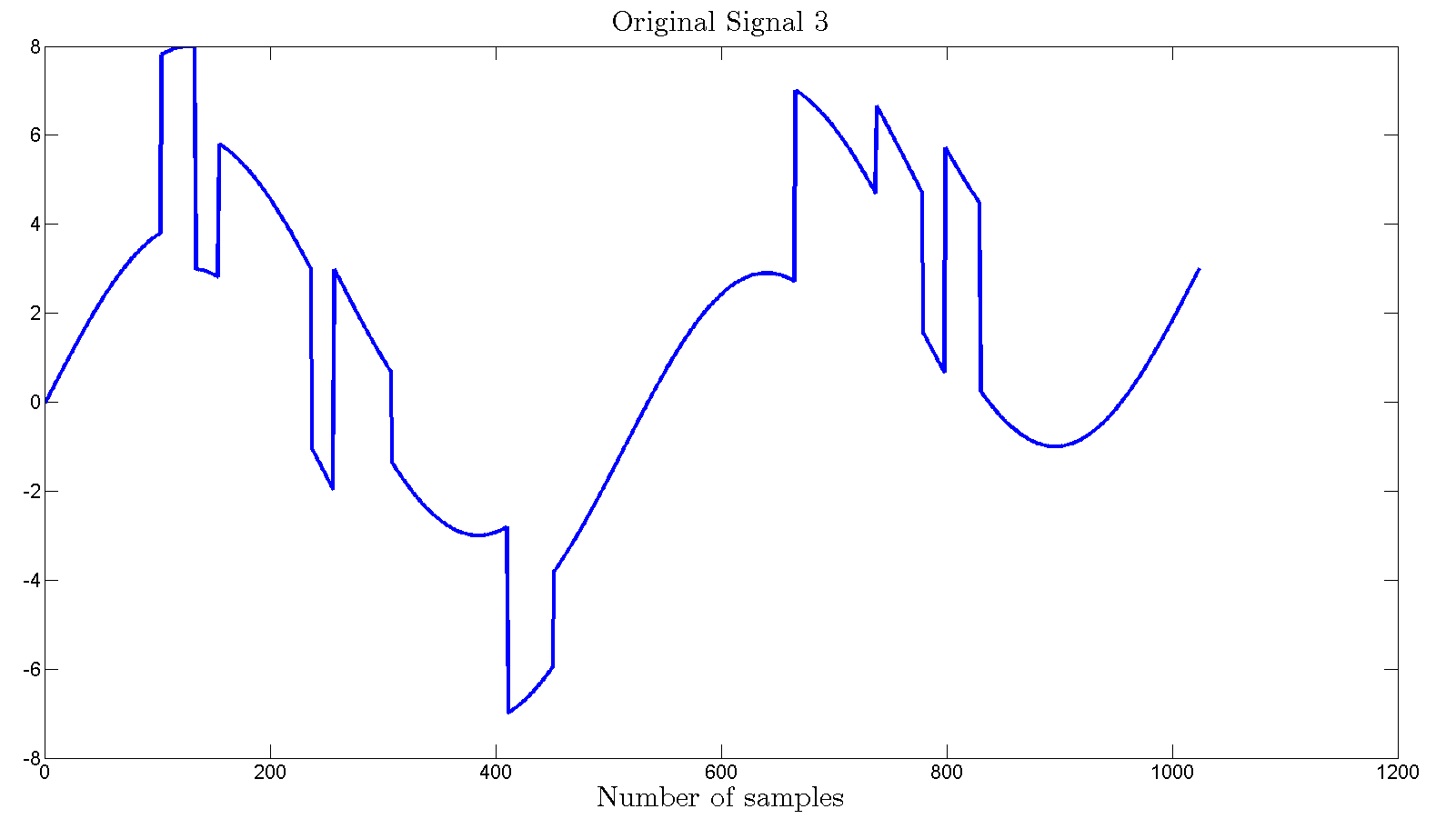}}
\subfigure[Signal 2]
{\label{Signal 4}\includegraphics[width=75mm]{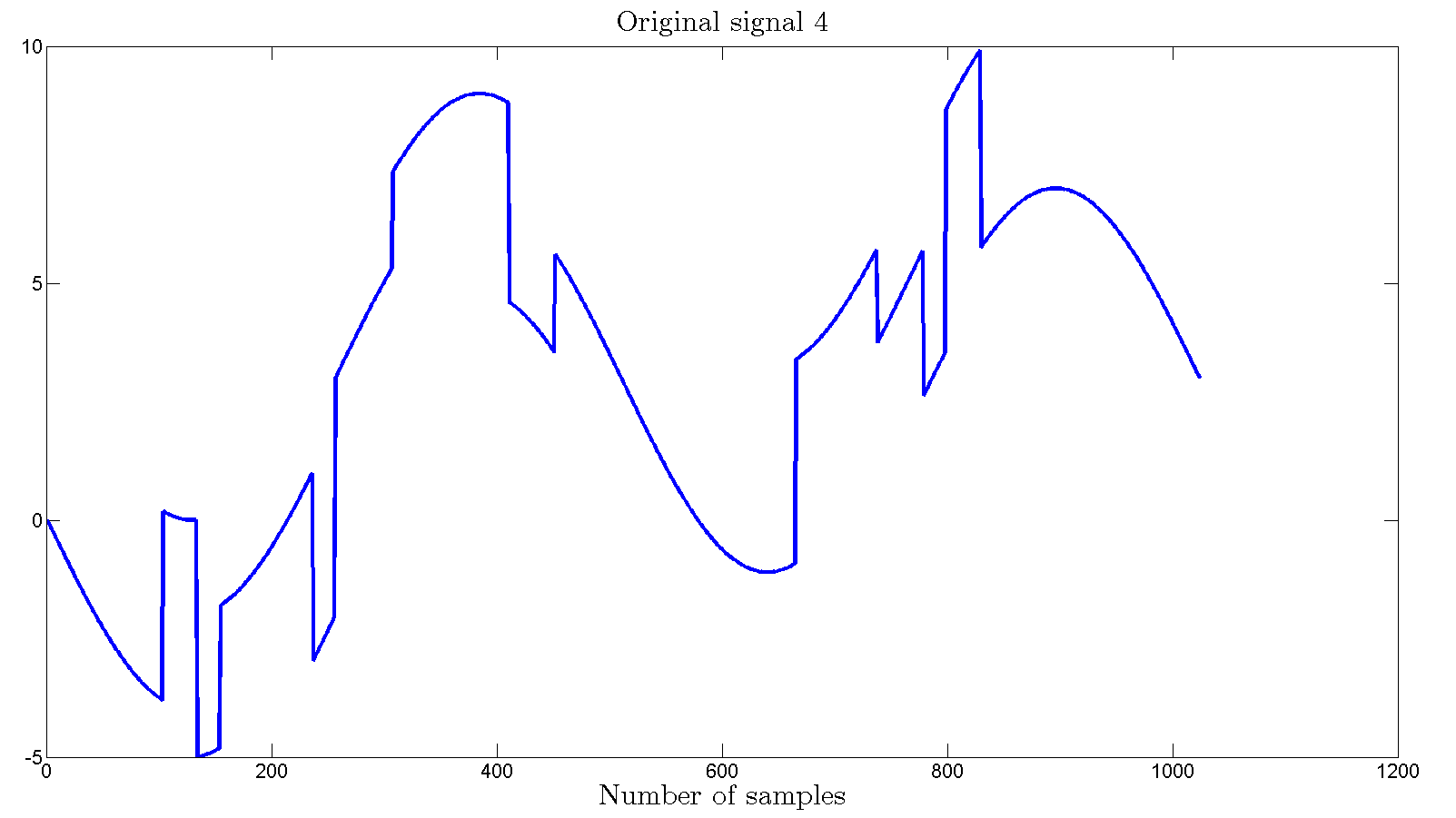}}
\subfigure[Blocks]
{\label{Signal 1}\includegraphics[width=75mm]{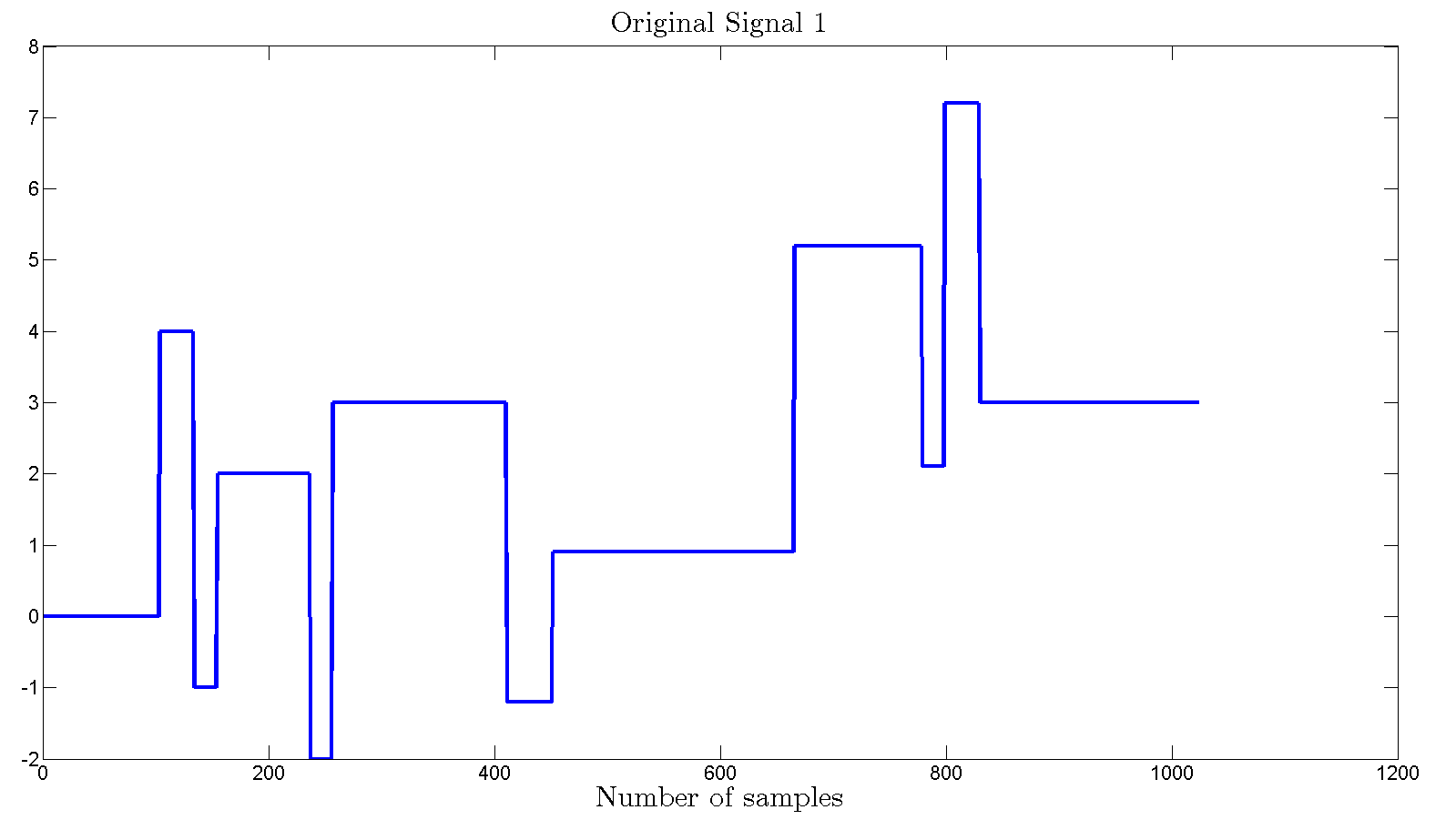}}
\subfigure[Heavy sine]
{\label{Signal 2}\includegraphics[width=75mm]{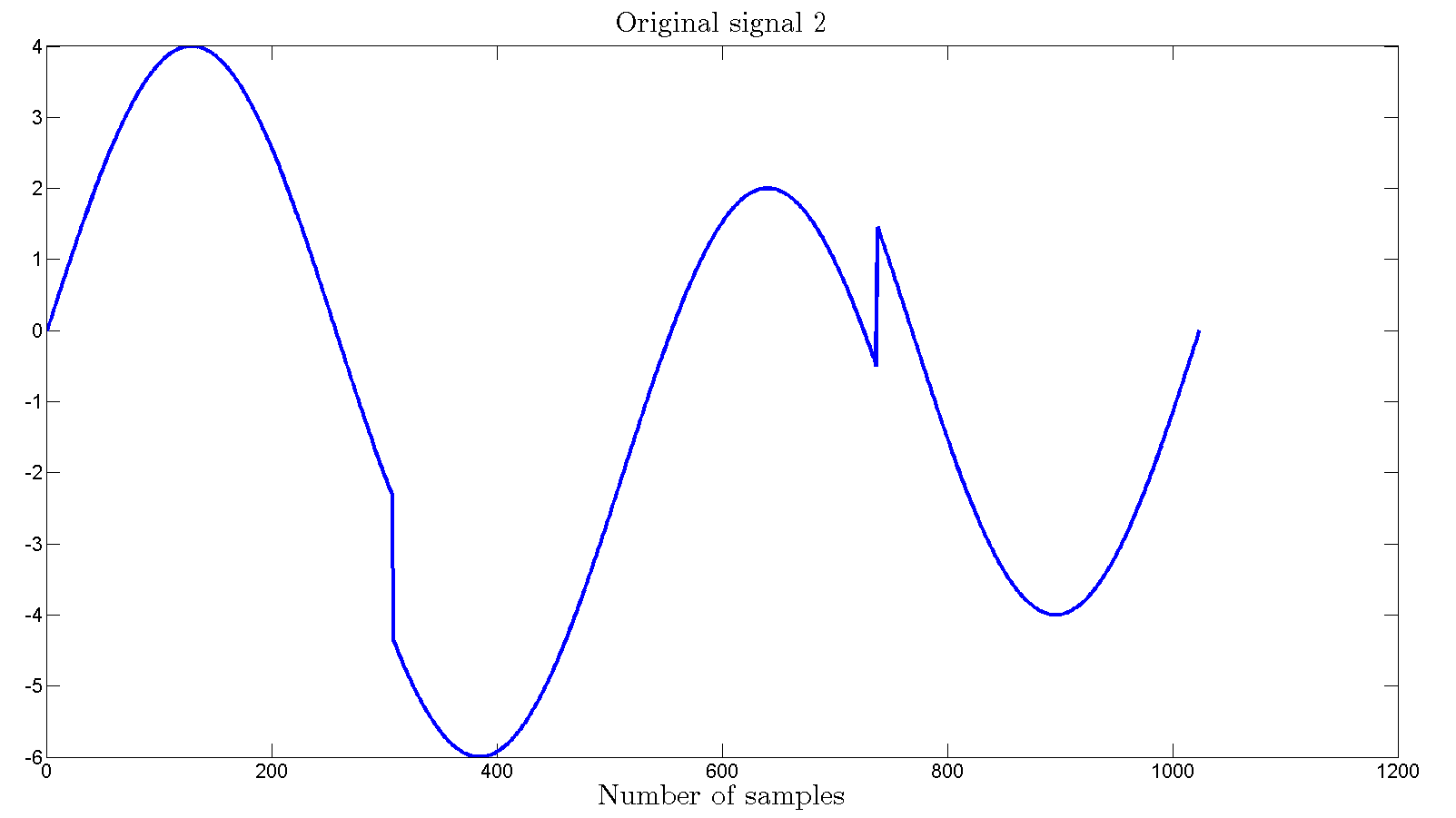}}
\subfigure[Piece-regular]
{\label{Piece-regular}\includegraphics[width=75mm]{piece.png}}
\subfigure[CUSP]
{\label{CUSP}\includegraphics[width=75mm]{CUSP.png}}
\caption{Signals which are used in the simulations.}
\label{app:OrigSig}
\end{figure}

\begin{table*}[ht!]
\begin{center}
\caption{Comparison of the results for denoising algorithms with Gaussian noise with $\sigma =$ 10, 20, and 30 \% of maximum amplitude of original signal.}
\label{tab:1}
\scalebox{0.9} {
\begin{tabular}{|c|c|c|c|c|c|c|c|}
\hline
\parbox[t]{1cm}{Signal}&\parbox[t]{1.5cm}{Input~SNR (dB)}&\parbox[t]{1.5cm}{\textbf{PES-$\ell_1$ Pyramid}}&\parbox[t]{1.5cm}{PES-$\ell_1$ Wavelet }&\parbox[t]{1.5cm}{MATLAB \cite{Mina,Rousseeuw}}&\parbox[t]{2.1cm}{Soft-threshold 3$\hat{\sigma}$}&\parbox[t]{2cm}{MATLAB ``rigrsure" \cite{Johnstone}}&\parbox[t]{2.3cm}{MATLAB ``mimimaxi" \cite{Donoho}}\\\hline\hline
Blocks&12.30&17.27&17.08&15.73&17.64&18.32&16.59\\\hline
Heavy sine&17.77&26.17&26.62&26.87&26.22&27.82&27.75\\\hline
Signal 1&13.09&18.43&18.10&16.63&17.80&19.18&17.41\\\hline
Signal 2&13.83&20.37&19.94&18.39&18.92&20.53&19.08\\\hline
Piece-Regular&12.32&18.53&18.05&16.41&18.16&19.35&17.66\\\hline
CUSP&16.29&32.58&29.40&30.43&28.72&29.10&29.99\\\hline\hline
Blocks&6.28&14.34&13.98&12.92&12.87&14.18&13.43\\\hline
Heavy sine&11.75&23.84&23.79&23.52&21.38&23.06&23.71\\\hline
Signal 1&7.07&15.70&15.28&14.00&13.30&15.15&14.42\\\hline
Signal 2&7.80&17.13&17.07&15.84&14.56&16.65&16.20\\\hline
Piece-Regular&6.27&15.24&14.47&13.11&13.14&14.94&13.99\\\hline
CUSP&10.25&28.24&24.89&25.04&23.27&23.48&24.44\\\hline\hline
Blocks&2.76&12.52&12.55&11.37&10.13&12.05&11.64\\\hline
Heavy sine&9.20&20.89&21.78&21.32&18.79&20.17&21.05\\\hline
Signal 1&3.56&13.50&13.65&12.37&10.44&13.06&12.72\\\hline
Signal 2&4.26&15.14&14.25&14.06&12.05&14.37&14.30\\\hline
Piece-Regular&2.77&13.21&12.70&11.37&9.94&12.43&12.05\\\hline
CUSP&6.73&25.10&23.47&21.73&19.67&19.69&21.02\\\hline\hline
Average&9.13&\textbf{19.68}&18.73&17.84&17.06&18.53&18.19\\\hline
\end{tabular}}
\end{center}
\end{table*}

\clearpage
\section*{Matlab Code}
\subsection*{PES-$\ell_1$ with pyramid method}
In the codes for PES-$\ell_1$ with pyramid method, first it starts with loading the original signals as:
\begin{lstlisting}[firstline=6,lastline=9]
x_orig = zeros(1024,6);
load ex4mwden
x_orig(:, 5) = load_signal('Piece-Regular', 1024);
x_orig(:, 6) = load_signal('cusp', 1024);
\end{lstlisting}
Then the white Gaussian noise is added as below:
\begin{lstlisting}[firstline=16,lastline=20]
    amp_perc = 0.1;
    for m = 1:kk
        sigma(m) = amp_perc*max(x_orig(:, m));
        NoisySignal(:, m) = x_orig(:, m) + sigma(m)*randn(size(x_orig(:, m)));
    end
\end{lstlisting}
which the noise standard deviation is determined with ``amp\_perc" which in our software it is 0.1, 0.2, and 0.3. Then the iteration number is determined according to the noise power. After all the signal will enter the denoising algorithm which is applied to the noisy signal in ``PES\_L1\_Pyramid.m" function, therefore:
\begin{lstlisting}[firstline=31,lastline=31]
    DenoisedSignal = PES_L1_Pyramid(iter, NoisySignal, kk);    
\end{lstlisting}
which ``iter" is the number of the iterations, ``NoisySignal" is the corrupted signal, and ``kk" is the number of the signals, which here we have six signals in our simulations.  In the main function of PES-$\ell_1$ denoising ``PES\_L1\_Pyramid.m", first the signal is passed through the high-pass filter and signal's high and low frequencies are separated, and then the PES-$\ell_1$ algorithm is applied to the high-passed signal. After that, the denoised high-pass signal is added to the unchanged low-pass signal and the main denoised signal is obtained. AS mentioned above, the performance of the algorithms are evaluated by SNR. which are calculated as:
\begin{lstlisting}[firstline=33,lastline=36]
    for d = 1:kk
        SNR_in(d) = snr(x_orig(:,d), NoisySignal(:,d));
        SNR_out(d) = snr(x_orig(:,d), DenoisedSignal(:,d));
    end
\end{lstlisting}
Since the additive noise is random then we have to run the codes repeatedly for enough times and average them to get the rational SNR value. Which averaging is done as:
\begin{lstlisting}[firstline=56,lastline=57]
SNR_In_Ave = mean(SNR_in1, 1)
SNR_Out_Ave = mean(SNR_out1, 1)
\end{lstlisting}
Then all the signals (original, noisy, and denoised) are ploted as:
\begin{lstlisting}[firstline=43,lastline=53]
kp = 0;
figure(2)
for f = 1:kk
    subplot(kk,3,kp+1), plot(x_orig(:,f)); axis tight;
    title(['Original signal ',num2str(f)])
    subplot(kk,3,kp+2), plot(NoisySignal(:,f)); axis tight;
    title(['Observed signal ',num2str(f)])
    subplot(kk,3,kp+3), plot(DenoisedSignal(:,f)); axis tight;
    title(['Denoised signal ',num2str(f)])
    kp = kp + 3;
end
\end{lstlisting}
The resulting plot are:
\begin{figure}[!htb]
\centering
\includegraphics[scale=0.35]{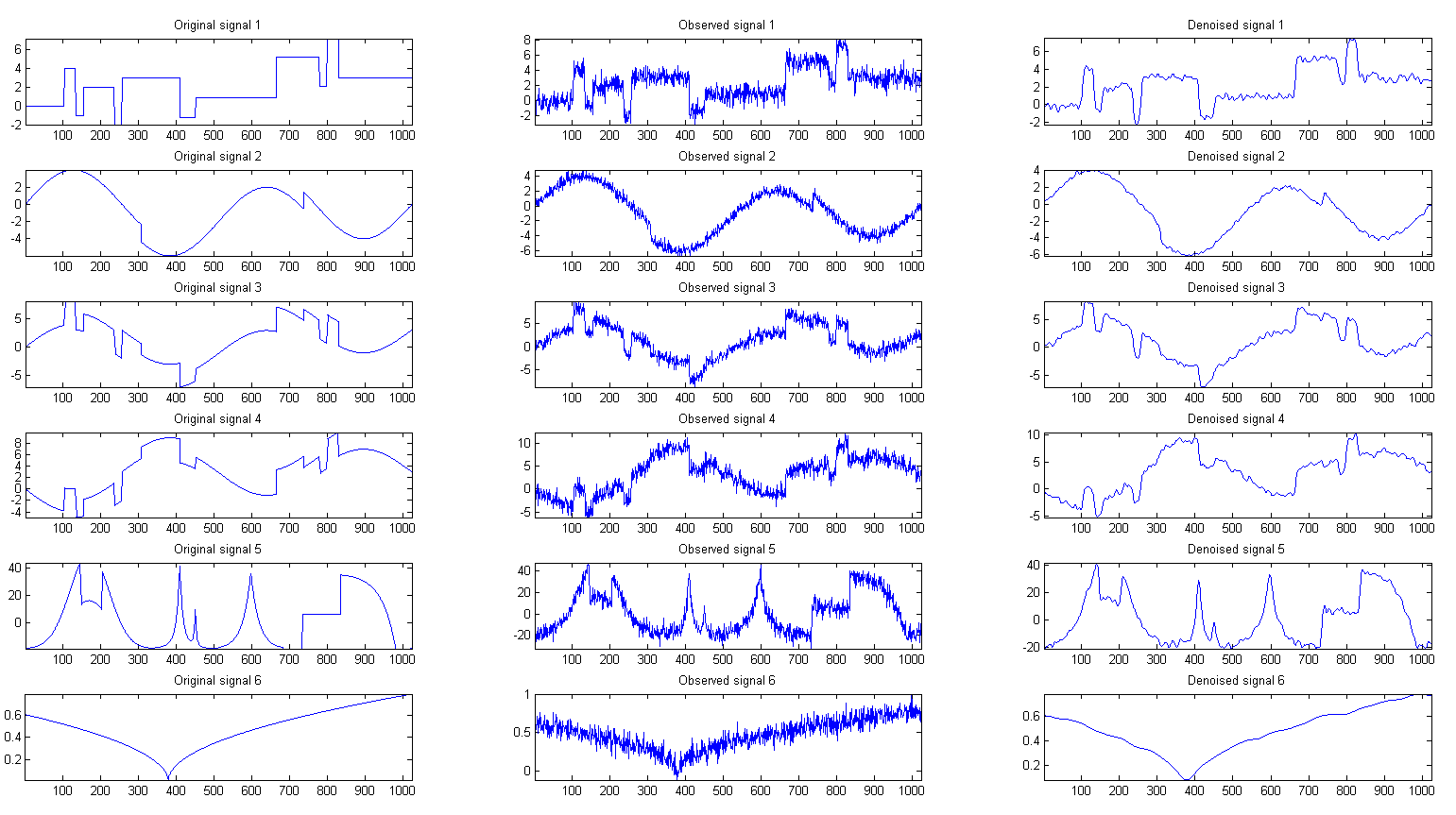}
\caption{All the original, noisy, and denoised signals for PES-$\ell_1$ with pyramid method.}
\label{sema1}
\end{figure}
and the SNRs are:
\begin{lstlisting}[firstline=75,lastline=87]
SNR_In_Ave =

    6.2994   11.7275    7.0773    7.7911    6.2919   10.2356


SNR_Out_Ave =

   14.3881   23.8448   15.6371   17.1841   15.3298   28.1067


ans =

   19.0818
\end{lstlisting}

\subsection*{PES-$\ell_1$ with wavelet method}
All the preliminary steps for this codes are as the same for PES-$\ell_1$ with pyramid method. Here, instead of ``PES\_L1\_Pyramid.m", the function for PES-$\ell_1$ with wavelet method ``PES\_L1\_Wavelet.m" is used. In this function the ``farras" filter bank is used for wavelet decomposition and the decomposition level is determined as explained in previous sections. It is done as:
\begin{lstlisting}[firstline=5,lastline=8]
        x_dwt1 =  circshift(NoisySignal(:, k), f-1); % Shifting the signal
        [af, sf] = farras; % Wavelet transforms filters
        J = Jk(k); % Decomposition level
        x_dwt = dwt_C(x_dwt1,J,af); % Wavelet decomposition
\end{lstlisting}
then the high subsignals are denoised with PES-$\ell_1$ algorithm and the low subsignal is transfered to the output without without any change, then the main denoised signal is reconstructed as follows;
\begin{lstlisting}[firstline=39,lastline=42]
        % Reconstructing the signal from its subsignals
        im_den(J+1) = x_dwt(J+1);
        x_idwt1 = idwt_C(im_den, J, sf);
        x_idwt2(:, f) = circshift(x_idwt1, -f+1); % Shift back
\end{lstlisting}
again the SNR is calculated as before, and the signals are plotted as follows:
\begin{lstlisting}[firstline=63,lastline=73]
kp = 0;
figure
for f = 1:kk
    subplot(kk,3,kp+1), plot(x_orig(:,f)); axis tight;
    title(['Original signal ',num2str(f)])
    subplot(kk,3,kp+2), plot(NoisySignal(:,f)); axis tight;
    title(['Observed signal ',num2str(f)])
    subplot(kk,3,kp+3), plot(DenoisedSignal(:,f)); axis tight;
    title(['Denoised signal ',num2str(f)])
    kp = kp + 3;
end
\end{lstlisting}
and the SNR values are averaged as before and the signals are plotted as following figure:
\begin{figure}[!htb]
\centering
\includegraphics[scale=0.35]{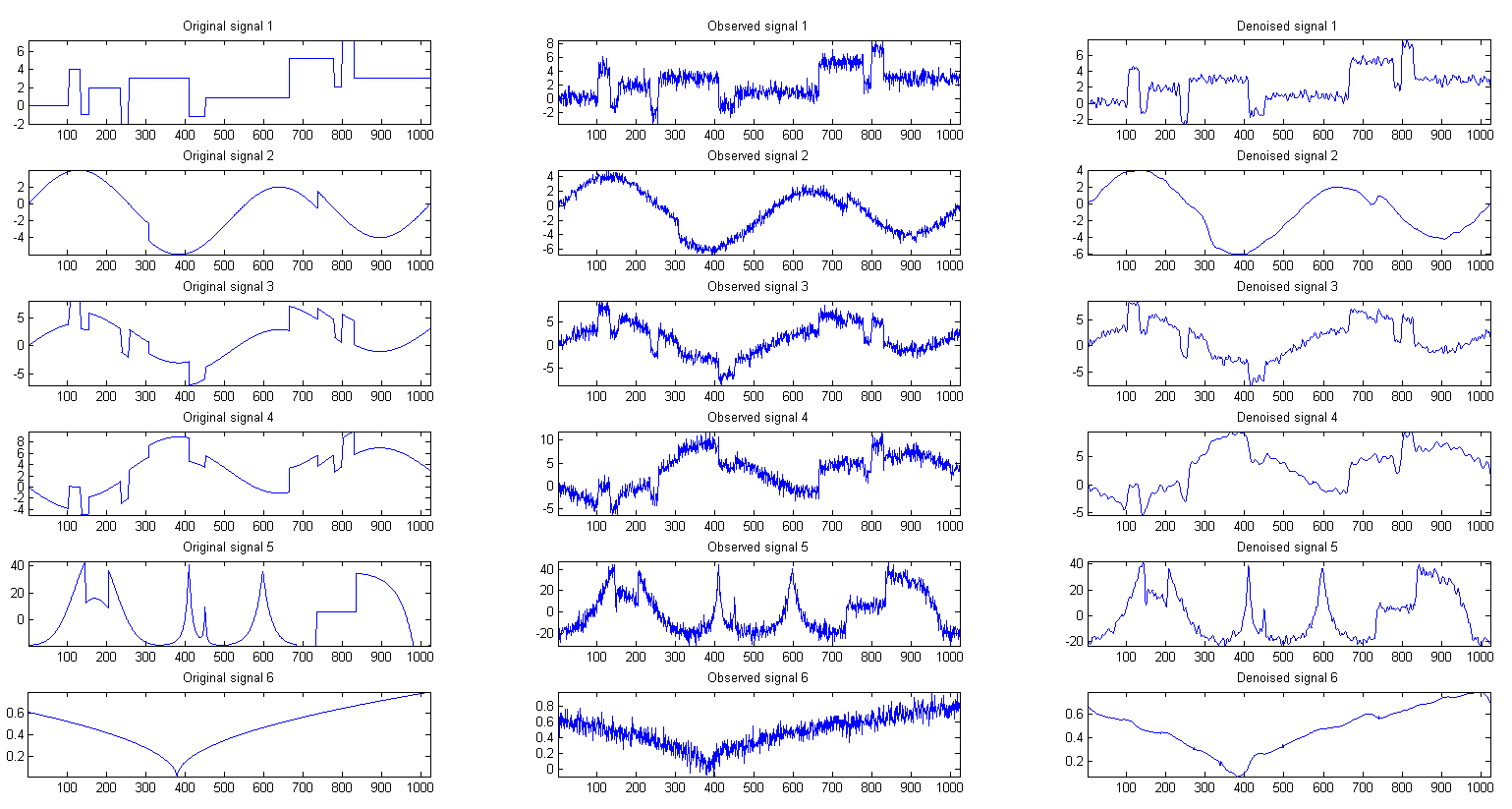}
\caption{All the original, noisy, and denoised signals for PES-$\ell_1$ with pyramid method.}
\label{sema1}
\end{figure}
and the SNRs are:
\begin{lstlisting}[firstline=89,lastline=101]
   SNR_In_Ave =

    6.2991   11.7385    7.0891    7.7825    6.2741   10.2704


SNR_Out_Ave =

   13.9855   23.7737   15.2842   17.0870   14.4843   24.8300


ans =

   18.2408
\end{lstlisting}

\clearpage
\bibliographystyle{IEEEtran}
\bibliography{Report}

\end{document}